\documentclass{amsart}

\usepackage{amscd,amsmath,xypic,amssymb,combelow,tikz-cd,etoolbox,calligra,mathrsfs,enumitem,mathtools,epigraph}

\usepackage{graphicx, amsfonts,amsthm, amscd,amsbsy,amstext, tikz, tensor}  
\usepackage{latexsym,ifthen,stmaryrd,fancyhdr,empheq,verbatim}
\mathtoolsset{showonlyrefs,showmanualtags}

\usepackage[russian,english]{babel}

\usepackage{hyperref}
\hypersetup{anchorcolor=black, linktocpage=true,
	colorlinks=true,
	citecolor=black,
	linkcolor=black,
	urlcolor=black,
	breaklinks=true,
	plainpages=true
}

\usetikzlibrary{positioning,shapes,shadows,arrows}

\usepackage[utf8x]{inputenc}

\usepackage[toc,page]{appendix}
\usepackage[mathscr]{euscript} %with mathscr option
\allowdisplaybreaks

\usepackage[russian,english]{babel}

\makeatletter
\patchcmd{\@settitle}{\uppercasenonmath\@title}{}{}{}
\makeatother

\xyoption{all}
\CompileMatrices

\makeatletter
\@addtoreset{equation}{section}
\makeatother

\newtheorem{theorem}[subsection]{Theorem}

\newtheorem{proposition}[subsection]{Proposition}
\newtheorem{lemma}[subsection]{Lemma}
\newtheorem{corollary}[subsection]{Corollary}

\newtheorem{definition}[subsection]{Definition}

\newtheorem{claim}[subsection]{Claim}
\newtheorem{example}[subsection]{Example}
\newtheorem{remark}[subsection]{Remark}

\def\loccitt{\emph{loc. cit.}}
\def\loccit{\emph{loc. cit. }}

\def\fsl{{\mathfrak{sl}}}

\def\BC{{\mathbb{C}}}
\def\BK{{\mathbb{K}}}

\def\BN{{\mathbb{N}}}
\def\BF{{\mathbb{F}}}

\def\BQ{{\mathbb{Q}}}
\def\BT{{\mathbb{T}}}
\def\BZ{{\mathbb{Z}}}

\def\CS{{\mathcal{S}}}

\def\CV{{\mathcal{V}}}
\def\CW{{\mathcal{W}}}

\def\vs{\varsigma}

\def\pt{\textrm{pt}}
\def\and{\textrm{ }\&\textrm{ }}

\def\esym{\emph{sym}}

\def\te{{e}}

\def\nn{{{\BN}}^I}

\def\UU{\mathbf{U}}
\def\UUp{\mathbf{U}^+}
\def\UUm{\mathbf{U}^-}

\def\tUU{{\widetilde{\mathbf{U}}}}
\def\tUUp{{\widetilde{\mathbf{U}}^+}}

\def\bs{{\boldsymbol{\vs}}}

\def\bn{\boldsymbol{n}}

\def\hdeg{\text{hdeg }}
\def\vdeg{\text{vdeg }}

\def\op{\text{op}}

\def\oij{\overrightarrow{ij}}
\def\oji{\overrightarrow{ji}}

\def\edge{e}
\def\tzeta{\widetilde{\zeta}}

\def\K{K_Q}
\def\Kloc{K_{Q,\text{loc}}}

\def\bE{\overline{E}}

\def\finite{\mathbf{k}}

\def\LHom{\text{LHom}}

\begin{document}

\title[Shuffle algebras for quivers as quantum groups]{\Large{\textbf{Shuffle algebras for quivers as quantum groups}}}

\author[Andrei Negu\cb t]{Andrei Negu\cb t}
\address[Andrei Negu\cb t]{MIT, Department of Mathematics, Cambridge, MA, USA}
\address{Simion Stoilow Institute of Mathematics, Bucharest, Romania}
\email{andrei.negut@gmail.com}

\author[Francesco Sala]{Francesco Sala}
\address[Francesco Sala]{Università di Pisa, Dipartimento di Matematica, Largo Bruno Pontecorvo 5, 56127 Pisa (PI), Italy}
\address{Kavli IPMU (WPI), UTIAS, The University of Tokyo, Kashiwa, Chiba 277-8583, Japan}
\curraddr{}
\email{\href{mailto:francesco.sala@unipi.it}{francesco.sala@unipi.it}}

\author[Olivier Schiffmann]{Olivier Schiffmann}
\address[Olivier Schiffmann]{Laboratoire de Mathématiques d'Orsay, Université de Paris-Saclay, Orsay, France}
\email{olivier.schiffmann@universite-paris-saclay.fr}

\maketitle

\begin{abstract} We define a quantum loop group $\UUp_Q$ associated to an arbitrary quiver $Q=(I,E)$ and maximal set of deformation parameters, with generators indexed by $I \times \BZ$ and some explicit quadratic and cubic relations. We prove that $\UUp_Q$ is isomorphic to the (generic, small) shuffle algebra associated to the quiver $Q$ and hence, by \cite{N}, to the localized $K$-theoretic Hall algebra of $Q$. For the quiver with one vertex and $g$ loops, this yields a presentation of the spherical Hall algebra of a (generic) smooth projective curve of genus $g$ (invoking the results of \cite{SV}). We extend the above results to the case of non-generic parameters satisfying a certain natural metric condition. As an application, we obtain a description by generators and relations of the subalgebra generated by absolutely cuspidal eigenforms of the Hall algebra of an arbitrary smooth projective curve (invoking the results of \cite{KSV}).

\end{abstract}

\setcounter{tocdepth}{1}
\tableofcontents

\section{Introduction}

\subsection{} Let $Q$ be a finite quiver, with vertex set $I$ and edge set $E$; edge loops and multiple edges are allowed. The Hall algebra of the category of representations of $Q$ over a finite field is well-known to contain a copy of the quantized enveloping algebra $U^+_q(\mathfrak{g}_Q)$, where $\mathfrak{g}_Q$ is the Kac-Moody Lie algebra associated to $Q$ (or the Bozec-Kac-Moody Lie algebra when $Q$ has edges loops). Cohomological Hall algebras associated to $Q$ for a Borel-Moore homology theory (including $K$-theory) were more recently introduced, in relation to Donaldson-Thomas theory on the one hand and Nakajima quiver varieties on the other hand (see \cite{KS, SVHilb, YZ}). More precisely, the \textit{$K$-theoretic Hall algebra} of $Q$ is the vector space:
$$
\K \coloneqq\bigoplus_{\bn \in \BN^I} K^{\BT}(\text{Rep}_{\bn} \, \Pi_Q)
$$
where $\Pi_Q$ is the preprojective algebra of $Q$, and $\text{Rep}_{\bn}\, \Pi_Q$ is the stack of complex $\bn$-dimensional representations of $\Pi_Q$. The vector space $\K$ is equipped with a natural Hall multiplication making it into an associative algebra \footnote{There is a specific choice of a line bundle involved in the definition of the multiplication; we refer to \cite{N} for details.}. Here $\BT$ is a torus acting in a Hamiltonian way on $\text{Rep}_{\bn} \;\Pi_Q$ by appropriately rescaling the maps attached to the arrows $e \in E$. The algebra $\K$ acts on the $\BT$-equivariant $K$-theory groups of Nakajima quiver varieties (see \cite{N2} for a review in our language):
$$
\mathcal{N}_{\boldsymbol{w}} = \bigsqcup_{\boldsymbol{v} \in \nn} \mathcal{N}_{\boldsymbol{v}, \boldsymbol{w}}
$$and is in fact the largest algebra thus acting via Hecke correspondences. When $Q$ is a finite type quiver, $\K$ is isomorphic to the positive half of the quantum loop algebra (in Drinfeld's sense) of $\mathfrak{g}_Q$. The situation of an affine quiver, in which case $\K$ is isomorphic to a quantum toroidal algebra, is studied in detail for the Jordan quiver in \cite{SVHilb}, for cyclic quivers in \cite{Negutthesis} and for arbitrary affine quivers in \cite{VV}. More generally, for a quiver without edge loops and a specific one-dimensional torus $\BT$, there is an algebra homomorphism:
\begin{equation}
\label{eqn:quantumafftokha}
U_q^+(L\mathfrak{g}_Q) \longrightarrow \K
\end{equation}
which recovers Nakajima's construction of representations of quantum affinizations of Kac-Moody algebras on the equivariant $K$-theory of quiver varieties. The map \eqref{eqn:quantumafftokha} is surjective (under some mild conditions on the torus action on $\text{Rep}_{\bn} \, \Pi_Q$), but it is not known to be injective in general (see \cite{VV}). Beyond these cases, however, very little is known. Moreover, even though $\text{Rep}_{\bn} \, \Pi_Q$ is equivariantly formal for any $\BT$ (and hence $K^{\BT}(\text{Rep}_{\bn} \, \Pi_Q)$ is free as a $K^{\BT}(\pt)$-module), the structure of $\K$ as an algebra depends in a rather subtle way on $\BT$. Note that there is a natural gauge action of the group $(\BC^*)^I$ on $\BT$, but as soon as $Q$ contains edge loops or multiple edges, the quotient of $\BT$ by this gauge group is nontrivial. 

\medskip

\subsection{} In the present paper, we consider the case when the torus $\BT = (\BC^*)^{|E|} \times \BC^*$ is as large as possible (each of the first $|E|$ copies of $\BC^*$ scale anti-diagonally the two coordinates of $\Pi_Q$ corresponding to a given edge, while the last copy of $\BC^*$ scales diagonally one half of the coordinates of $\Pi_Q$) and we work over the fraction field: 
\begin{equation}
\label{eqn:ground field}
\mathbb{F} = \text{Frac}(K^{\BT}(\pt))= \mathbb{Q}(q,t_e)_{e \in E}
\end{equation}
Our main result provides an explicit description of: 
$$
\Kloc\coloneqq \K \bigotimes_{K^{\BT}(\pt)} \mathbb{F}
$$
by generators and relations which we will now summarize. Let $\overline{E}$ be the ``double" of the edge set $E$, i.e. there are two edges $e = \oij$ and $e^* = \oji$ in $\bE$ for every edge $e = \oij \in E$. The set $\overline{E}$ is equipped with a canonical involution $e \leftrightarrow e^*$. We extend the notation $t_e$ to an arbitrary $e \in \overline{E}$ by the formula: 
\begin{equation}
\label{eqn:inverse t}
t_{e^*}= \frac q{t_e}
\end{equation}
for any $e \in E$. For any $i,j \in I$, consider the rational function \footnote{Note that this is actually the rational function denoted by $\zeta_{ij}'$ in \cite{N}.}:
\begin{equation}
\label{eqn:def zeta}
\zeta_{ij}(x) = \left(\frac {1-xq^{-1}}{1-x} \right)^{\delta_j^i} \prod_{e = \oij \in E} \left(\frac 1{t_e} - x \right) \prod_{e = \oji \in E} \left(1 - \frac {t_e}{qx} \right)
\end{equation}
and set:
$$
\tzeta_{ij}(x) = \zeta_{ij}(x) \cdot (1-x)^{\delta_j^i}
$$
Let $\UUp_Q$ be the $\BF$-algebra generated by elements $e_{i,d}$ for $i \in I, d \in \BZ$ subject to the following set of quadratic and cubic relations, in which we set 
$$
e_i(z) = \sum_{d \in \BZ} \frac {\te_{i,d}}{z^{d}}
$$

\medskip

\begin{itemize}[leftmargin=*]

\item For any pair $(i,j) \in I^2$, the \textbf{quadratic} relation:
\begin{equation}
\label{eqn:rel quadintro}
e_i(z) e_j (w) \tzeta_{ji} \left( \frac wz \right) z^{\delta_j^i} = e_j(w) e_i(z) \tzeta_{ij} \left(\frac zw \right)  (-w)^{\delta_j^i}
\end{equation}

\item For any edge $\overline{E} \ni e = \oij$, the \textbf{cubic} relation:
\begin{multline}
\label{eqn:serre a}
\frac {\tzeta_{ii}\left(\frac {x_2}{x_1} \right) \tzeta_{ji} \left( \frac {y}{x_1} \right) \tzeta_{ji} \left( \frac {y}{x_2} \right)}{\left(1- \frac {x_2}{x_1q}\right) \left(1-\frac {yq}{x_2t_{\edge}} \right)} \cdot e_i(x_1) e_i(x_2) e_j(y)  \\
+ \frac {\tzeta_{ii}\left(\frac {x_1}{x_2} \right) \tzeta_{ji} \left( \frac y{x_2} \right) \tzeta_{ij} \left( \frac {x_1}y \right) \left(- \frac {x_2 t_{\edge}}{y} \right) \left(- \frac y{x_1} \right)^{\delta_j^i}}{\left(1- \frac {yq}{x_2t_{\edge}}\right) \left(1-\frac {x_1 t_{\edge}}{y} \right)} \cdot e_i(x_2) e_j(y) e_i(x_1)  \\
+ \frac {\tzeta_{ii}\left(\frac {x_2}{x_1} \right) \tzeta_{ij} \left( \frac {x_1}y \right) \tzeta_{ij} \left( \frac {x_2}y \right) \left( \frac {x_2t_{\edge}}{yq} \right)\left( \frac {y^2}{x_1x_2} \right)^{\delta_j^i}}{\left(1- \frac {x_2}{x_1 q}\right) \left(1-\frac {x_1 t_{\edge}}{y} \right)} \cdot e_j(y) e_i(x_1) e_i(x_2) = 0
\end{multline}

\end{itemize}

\medskip

\begin{theorem}
\label{thm:KHA} 

There is an algebra isomorphism $K_{Q,\emph{loc}} \simeq \mathbf{U}^+_Q$. \medskip

\end{theorem}

\noindent When $Q$ is a tree, the quotient of $ (\BC^*)^{|E|} \times \BC^*$ by the action of the gauge group $(\BC^*)^{|I|}$ is one-dimensional, hence up to renormalization $K_Q$ can be defined as $\BC^*$-equivariant $K$-theory (in other words, we do not lose any information by assuming that $t_e = q^{\frac 12}$, for all $e \in E$). In addition, one can check that in the case of an $A_2$ quiver, the cubic relations \eqref{eqn:serre a} are equivalent to the standard $q$-Serre relations (see Example~\ref{ex:a1}). With this in mind, Theorem \ref{thm:KHA} implies: \medskip

\begin{theorem} 
\label{thm:tree}

Suppose that $Q$ is a tree, and that $\BT$ scales the symplectic form on $\emph{Rep}_{\bn} \,\Pi_Q$ nontrivially. Then the localization: 
$$
U_q^+(L\mathfrak{g}_Q) \bigotimes_{K^{\BT}(\emph{pt})} \mathbb{F} \longrightarrow K_{Q,\emph{loc}}
$$
of the map \eqref{eqn:quantumafftokha} is an algebra isomorphism. \medskip

\end{theorem}

\noindent Cohomological Hall algebras of quivers are known (at least in the case of Borel-Moore homology and $K$-theory) to embed in a suitable big \textbf{shuffle algebra} $\CV_Q$, whose multiplication encodes the structure of $Q$ (see \cite{SVgen, VV, YZ}). In the $K$-theoretic case and for maximal $\BT$, recent work (\cite{N, Zhao}) identified the image of this embedding as the small shuffle algebra $\CS_Q \subset \CV_Q$ determined by the so-called 3-variable \textbf{wheel conditions}. These wheels conditions feature for instance in \cite{Enriquez} where they are derived from purey algebraic formal identities involving products of delta functions (see also \cite{Ding-Jing}) and more recently in \cite{NT}.
Theorem~\ref{thm:KHA} is a direct corollary of the following theorem, which is the main result of the present paper. \medskip

\begin{theorem}
\label{thm:intro}

There is an algebra isomorphism $\mathbf{U}_Q^+ \cong \CS_Q$. \medskip

\end{theorem}

\noindent For a general $Q$ and a general choice of $\BT$ (which satisfies some mild conditions), there is a chain of algebra homomorphisms:
$$
\mathbf{U}_Q^+ \longrightarrow \Kloc \longrightarrow \CS_Q
$$
The content of Theorem~\ref{thm:intro} is that these maps are all isomorphisms for $\BT$ maximal. \medskip

\noindent Our main tool to prove Theorem~\ref{thm:intro} is the combinatorics of words developed in \cite{N} (which was in turn influenced by \cite{NT}, and further back, by the seminal work of \cite{LR,L,R}). It would be interesting to extend the above result to the case of a smaller torus $\BT$; this would necessitate some more complicated wheel conditions, and result in higher degree relations in $\UUp_Q$ (see the last section of \cite{N}). For instance, when two vertices are joined by more than one edge, it is customary to set the corresponding weights of the torus action to be in a geometric progression, as this yields the $q$-Serre relations (which are of degree two more than the number of edges, so more complicated than cubic in the case of multiple edges). 

\medskip

\subsection{} Let us mention an important application of Theorem~\ref{thm:intro}. When $Q$ is the quiver $S_g$ with one vertex and $g$ loops, it is known by combining \cite{SV} and \cite{N} that the \textbf{spherical} subalgebra of the \textbf{Hall algebra} $ \mathbf{H}^{sph}_X$ of the category of coherent sheaves on a genus $g$ curve $X$ defined over the finite field $\finite$ of $q^{-1}$ elements is isomorphic to $\K$ (extended by a commutative Cartan subalgebra). To make the previous statement precise, the equivariant parameters $t_1, \ldots, t_g,q/t_1,\dots,q/t_g$ must be set equal to the inverses of the Weil numbers $\sigma_1, \ldots, \sigma_g, \overline{\sigma_1}, \dots, \overline{\sigma_g}$ of $X$. Thus, the rational function \eqref{eqn:def zeta} for the quiver $Q = S_g$ corresponds to the renormalized zeta function of the curve $X$:
\begin{equation}
\label{eqn:zeta g}
\zeta_X(x)=\frac{1-xq^{-1}}{1-x}\prod_{e=1}^g (\sigma_e-x)(1-\overline{\sigma_e}x^{-1})
\end{equation}
For any $e=1, \ldots, g$, set:
$$
Q_e(z_1,z_2,z_3)=\prod_{1 \leq i < j \leq 3}\prod_{f \neq e} \left(\sigma_f-\frac{z_j}{z_i}\right)\left(1-\overline{\sigma_f}\frac{z_i}{z_j}\right)
$$
Then we have the following result (see Section~\ref{sec:genusg} for details).

\medskip

\begin{theorem}\label{thm:genusg}

When the curve $X$ has distinct Weil numbers (we will refer to this as the ``generic" case), its spherical Hall algebra $\mathbf{H}^{sph}_X$ is generated by elements $\kappa^{\pm 1}$, $\theta_{0,l}, 1^{vec}_{d}$ for $l\geq 1, d \in \BZ$, subject to the following set of relations:
\begin{align}\label{eqn:genusg1}
H^+(z)H^+(w)&=H^+(w)H^+(z)\\
\label{eqn:genusg2}
E(z)H^+(w)&=H^+(w)E(z) \frac{\zeta_X\left(\frac{z}{w}\right)}{\zeta_X\left(\frac{w}{z}\right)} \qquad (\text{for\;}|w| \gg |z|)
\\
\label{eqn:genusg3}
E(z)E(w)\zeta_X\left(\frac{w}{z}\right) &= E(w)E(z)\zeta_X\left(\frac{z}{w}\right) 
\end{align}
and for all $e=1, \ldots, g$ and $m \in \BZ$ the relation:
\begin{equation}
\label{eqn:genusg4}
\Big[(xyz)^m(x+z)(xz-y^2)Q_e(x,y,z) E(x)E(y)E(z) \Big]_{\emph{ct}} = 0
\end{equation}
(see \eqref{eqn:constant term} for the notation $[\dots]_{\emph{ct}}$). In the formulas above, we set:
$$
E(z)=\sum_{d \in \BZ} 1^{vec}_{d}z^{-d}, \qquad H^+(z)=\kappa\left( 1 + \sum_{l \geq 1} \theta_{0,l}z^{-l}\right)
$$
\end{theorem}

\medskip

\noindent
For $g=0$, one recovers the defining relations for $U_q^+(L\mathfrak{sl}_2)$ (see \cite{K}), while for $g=1$ one gets the relations describing the elliptic Hall algebra of \cite{BS}. Note that there is a slight discrepancy between \eqref{eqn:serre a} and \eqref{eqn:genusg4}; this comes from the fact that we added in the Cartan loop generators $\{\theta_{0,l}\}_{l \in \mathbb{N}}$, see Section~\ref{sec:genusg}. From the point of view of function field automorphic forms, Theorem~\ref{thm:genusg} says that in the case of the function field of a generic curve, Eisenstein series for the group $GL(n)$ which are induced from the trivial character of the torus satisfy, in addition to the celebrated functional equation (which is equivalent to \eqref{eqn:genusg3}), $g$ families of cubic relations and \emph{no} higher degree relations.

\medskip

\noindent One might wonder what happens for a \textit{non-generic} curve, or for the entire Hall algebra $\mathbf{H}_X$ rather than the spherical subalgebra $\mathbf{H}^{sph}_X$. In Section \ref{sec:wholeHall}, we answer both of these questions using a version of Theorem~\ref{thm:intro} that holds for equivariant parameters $t_e$ which satisfy a certain metric condition (\textbf{Assumption \begin{otherlanguage*}{russian}Ъ\end{otherlanguage*}} introduced in \cite{N}). As it turns out, this metric condition applies to the inverse Weil numbers of $X$ due to the (generalized) Riemann hypothesis and the functional equation for the zeta function (and more generally for the Rankin-Selberg L-functions attached to a pair of absolutely cuspidal eigenforms); this also allows us to give, using \cite{KSV}, a complete presentation for an arbitrary curve of the subalgebra $\mathbf{H}^{abs}_X$ generated by the coefficients of all absolutely cuspidal eigenforms (Corollary \ref{cor:abscusp}). As we show, the structure of this algebra only depends on the various orders of vanishing of the Rankin-Selberg L-functions attached to pairs of absolutely cuspidal eigenforms.

\medskip

\subsection{} The plan of the present paper is the following. From now on, we will fix the quiver $Q$ and write simply $\UUp$ instead of $\UUp_Q$. \medskip

\begin{itemize}[leftmargin=*]

\item In Section~\ref{sec:big}, we consider the $\mathbb{F}$-algebra $\tUUp$ generated by elements $\{\te_{i,d}\}_{i \in I, d \in \BZ}$ modulo the quadratic relations \eqref{eqn:rel quadintro} and check that it is naturally dual to the so-called big shuffle algebra $\CV$, see Proposition~\ref{prop:pairing}:
\begin{equation}
\label{eqn:pairing intro 1}
\tUUp \otimes \CV^{\op} \xrightarrow{\langle \cdot, \cdot \rangle} \BF.
\end{equation}

\item In Section~\ref{sec:small}, we recall the shuffle algebra $\CS \subset \CV$ of \cite{N}, which is cut out by the 3-variable wheel conditions \eqref{eqn:wheel}. We show that these wheel conditions arise by pairing with certain cubic elements that will be defined in Proposition~\ref{prop:wheel pairing}:
\begin{equation}
\label{eqn:elements intro}
\Big\{ A^{(\edge)}_d \Big\}_{\edge \in \overline{E}, d \in \BZ} \in \tUUp
\end{equation}
This allows us to prove that \eqref{eqn:pairing intro 1} descends to a non-degenerate pairing:
\begin{equation}
\label{eqn:pairing intro 2}
\UUp \otimes \CS^{\op} \xrightarrow{\langle \cdot, \cdot \rangle} \BF
\end{equation}
where $\UUp$ is the quotient of $\tUUp$ by the ideal generated by the elements \eqref{eqn:elements intro}. Comparing \eqref{eqn:pairing intro 2} with the pairing:
$$
\CS \otimes \CS^{\op} \xrightarrow{\langle \cdot, \cdot \rangle} \BF
$$
that was studied in \cite{N} allows us to conclude that $\UUp \cong \CS$, thus establishing Theorem \ref{thm:intro}. \medskip

\item In Section~\ref{sec:double}, we provide a definition of a natural Drinfeld-type double $\UU$ of $\UUp$.\medskip

\item In Section~\ref{sec:specialize}, we consider specializations of the shuffle algebra $\CS$ when the equivariant parameters satisfy Assumption \begin{otherlanguage*}{russian}Ъ\end{otherlanguage*} of Definition \ref{def:assumption}. We extend the main results (i.e. the presentation by generators and relations) in this context; this involves some new families of (still cubic) relations. \medskip

\item In Section~\ref{sec:genusg}, we recall the basic notions concerning the spherical Hall algebra $\mathbf{H}_X^{sph}$ of a generic genus $g$ smooth projective curve $X/\finite$ and prove Theorem~\ref{thm:genusg}. \medskip

\item Finally, in Section~\ref{sec:wholeHall}, we use Section~\ref{sec:specialize} to extend the results of Section~\ref{sec:genusg} to an arbitrary curve $X$, and to the (much) larger subalgebra $\mathbf{H}^{abs}_X \subset \mathbf{H}_X$ (corresponding by Langlands duality to all \textit{geometrically} irreducible local systems). 

\end{itemize}

\bigskip

\centerline{\textbf{Acknowledgements}}

\medskip

\noindent We would like to thank Alexander Tsymbaliuk for many interesting discussions and substantial feedback. A.N. gratefully acknowledges NSF grants DMS-$1760264$ and DMS-$1845034$, as well as support from the Alfred P.\ Sloan Foundation and the MIT Research Support Committee. The work of F.~S. was partially supported by JSPS KAKENHI Grant Number JP21K03197. He acknowledges the MIUR Excellence Department Project awarded to the Department of Mathematics, University of Pisa, CUP I57G22000700001. Moreover, he is a member of GNSAGA of INDAM.

Furthermore, the final revision of the paper took place under the MIT-UNIPI Project (XI call). 

\bigskip

\section{The big shuffle algebra and the quadratic quantum loop group} 
\label{sec:big}

\medskip

\subsection{} Let us introduce the big shuffle algebra $\CV$ and quadratic quantum loop group $\tUUp$ associated to a quiver. We will define a perfect pairing between them, as well as a homomorphism $\widetilde{\Upsilon}\colon \tUUp \to \CV$. 

\begin{remark}

The reader might keep in mind the following analogy: $\tUUp$ is like a Verma module, $\CV$ is like the corresponding dual Verma module, and $\widetilde{\Upsilon}$ is like the canonical map between them. As we will see in Section \ref{sec:small}, the image of $\widetilde{\Upsilon}$ will be the quantum loop group $\UUp$ that we are interested in, just like the image of the map from a Verma to the dual Verma is the irreducible highest weight representation. 

\medskip

\noindent Alternatively, the more geometry-oriented reader may also think about the inclusion of an affine open subset $j: O \hookrightarrow X$ of an algebraic variety $X$ together with a local system $\mathcal{L}$ on $O$: $\tUUp$ is like $j_!\mathcal{L}$, $\CV$ is like $j_*\mathcal{L}$, and $\UUp$ is like the intermediate extension $j_{!*}\mathcal{L}$, which is the image of the natural morphism $j_!\mathcal{L} \to j_*\mathcal{L}$.

\end{remark} \medskip

\noindent Throughout the present paper, we fix a finite quiver $Q$ with vertex set $I$ and edge set $E$. The notation $e = \oij$ will mean ``$e$ is an arrow going from $i$ to $j$". Let:
\begin{equation}
\label{eqn:number of edges} 
\#_{\oij} = |\text{arrows }\oij| \ , \qquad \#_{ij}= \#_{\oij} + \#_{\oji}
\end{equation}
We identify $\BZ^I$ with $\bigoplus_{i \in I} \BZ \bs_i$, where $\bs_i$ is the vector with a single 1 at the $i$-th spot, and zeroes everywhere else. In terms of quiver representations, $\bs_i$ is the dimension vector of the simple quiver representation $S_i$ supported at $i$. The set of natural numbers $\BN$ will be considered to include 0. \medskip

\subsection{} We begin by recalling the big shuffle algebra, defined over the field \eqref{eqn:ground field}: \medskip

\begin{definition}
\label{def:shuf big}

The \textbf{big shuffle algebra} is the vector space:
\begin{equation}
\label{eqn:big shuffle}
\CV = \bigoplus_{\bn \in \nn} \BF \left[\dots,z^{\pm 1}_{i1},\dots,z^{\pm 1}_{in_i},\dots \right]^{\esym}_{i \in I}
\end{equation}
endowed with the associative product:
\begin{equation}
\label{eqn:shuf prod}
R(\dots,z_{i1},\dots,z_{in_i},\dots) * R'(\dots,z_{i1},\dots,z_{in'_i},\dots) = 
\end{equation}
$$
\emph{Sym} \left[\frac {R(\dots,z_{i1},\dots,z_{in_i},\dots)R'(\dots,z_{i,n_i+1},\dots,z_{i,n_i+n_i'},\dots)}{\prod_{i\in I} n_i! \prod_{i\in I} n_i'!}\mathop{\prod^{i,j \in I}_{1 \leq a \leq n_i}}_{n_{j} < b \leq n_{j}+n'_{j}} \zeta_{ij} \left( \frac {z_{ia}}{z_{jb}} \right) \right]
$$
and the unit being the function 1 in zero variables ($\zeta_{ij}(x)$ is defined in \eqref{eqn:def zeta}). \medskip

\end{definition}

\noindent In Definition~\ref{def:shuf big}, ``sym" refers to the set of Laurent polynomials which are symmetric with respect to the variables $z_{i1},\dots,z_{in_i}$ for each $i \in I$ separately, and ``Sym" denotes symmetrization with respect to the variables $z_{i1},\dots,z_{i,n_i+n_i'}$ for each $i \in I$ separately. Note that even though the right-hand side of \eqref{eqn:shuf prod} seemingly has simple poles at $z_{ia} - z_{ib}$ for all $i \in I$ and all $a < b$, these poles vanish when taking the Sym, as the orders of such poles in a symmetric rational function must be even. \medskip

\noindent The big shuffle algebra is graded by $\nn \times \BZ$, where $\nn$ keeps track of the number of variables $\{n_i\}_{i \in I}$ of a Laurent polynomial $R$, while $\BZ$ keeps track of the homogeneous degree of $R$. Under these circumstances, we will write:
\begin{equation}
\label{eqn:deg}
\deg R = (\bn,d)
\end{equation}
and refer to:
\begin{equation}
\label{eqn:hdeg vdeg}
\hdeg R = \bn \qquad \text{and} \qquad \vdeg R = d
\end{equation}
as the horizontal and vertical degree of $R$, respectively. The graded pieces of the algebra $\CV$ will be denoted by $\CV_{\bn}$ (when grading by horizontal degree only) and by $\CV_{\bn,d}$ (when grading by both horizontal and vertical degree). \medskip

\subsection{} Let us now introduce the quadratic quantum loop group associated to $Q$. \medskip

\begin{definition}
\label{def:big quantum}

The \textbf{quadratic quantum loop group} $\tUUp$ is the $\BF$-algebra generated by symbols: 
$$
\{\te_{i,d}\}_{i \in I, d \in \BZ} 
$$
modulo the quadratic relations:
\begin{equation}
\label{eqn:rel quad}
\te_i(z) \te_j (w) \zeta_{ji} \left( \frac wz \right) = \te_j(w) \te_i(z) \zeta_{ij} \left(\frac zw \right) 
\end{equation}
for all $i,j \in I$, where $\te_i(z) = \sum_{d \in \BZ} \te_{i,d}z^{-d}$ . \medskip

\end{definition}

\noindent The meaning of relation \eqref{eqn:rel quad} is that one cancels denominators (which yields the equivalent relation \eqref{eqn:rel quadintro}) and then equates the coefficients of every $\{z^kw^l\}_{k,l \in \BZ}$ in the left and right-hand sides, thus yielding relations between the generators $\te_{i,d}$~:
\begin{multline}
\label{eqn:rel quad 2}
\te_{i,a} \te_{j,b} + \sum_{\bullet = -\#_{ij} - \delta_j^i}^{- 1} \text{coeff} \cdot \te_{i,a+\bullet} \te_{j,b-\bullet} = \\ 
= \gamma \cdot \te_{j,b+\#_{ij}+\delta_j^i} \te_{i,a-\#_{ij} -\delta_j^i} + \sum_{\bullet = 0}^{\#_{ij}+\delta_j^i - 1} \text{coeff} \cdot  \te_{j,b+\bullet} \te_{i,a-\bullet}
\end{multline}
where ``coeff" denotes various coefficients in $\BF$, arising from the numerator of the functions $\zeta_{ij}$ that appear in \eqref{eqn:rel quad}; observe that the constant $\gamma$ in the right-hand side of \eqref{eqn:rel quad 2} is different from 1 if $i = j$, which will be important in Proposition \ref{prop:proto}. Just like $\CV$, the algebra $\tUUp$ is also $\nn \times \BZ$-graded, with:
$$
\deg \te_{i,d} = (\bs_i,d).
$$
Therefore, we will use the terms ``horizontal degree" and ``vertical degree" pertaining to $\tUUp$ as well, in accordance with \eqref{eqn:hdeg vdeg}. \medskip

\begin{proposition}
\label{prop:big hom}

The assignment $\te_{i,d} \mapsto z_{i1}^d$, for $i \in I, d \in \BZ$ induces an algebra homomorphism:
\begin{equation}
\label{eqn:hom big shuffle}
\widetilde{\Upsilon}\colon \tUUp \longrightarrow \CV.
\end{equation}
\end{proposition}

\medskip

\noindent The proof of Proposition~\ref{prop:big hom} is a straightforward computation (indeed, all one needs to show is that relations \eqref{eqn:rel quad} hold in the big shuffle algebra) which we leave as an exercise to the interested reader. We remark that \eqref{eqn:hom big shuffle} is neither injective, nor surjective: its image is the small shuffle algebra, and its kernel is generated by cubic relations, as we will show in Section~\ref{sec:small}. \medskip

\subsection{}
As one of our main tools, we will next define a pairing between $\tUUp$ and $\CV$. Let $Dz = \frac {dz}{2\pi i z}$. Whenever we write $\int_{|z_1| \gg \dots \gg |z_n|}$ we are referring to a contour integral taken over concentric circles around the origin in the complex plane (i.e. an iterated residue at $0$, or $\infty$). The following result is close to \cite[Proposition~3.3]{N}. \medskip

\begin{proposition}
\label{prop:pairing}
There is a well-defined pairing:
\begin{equation}
\label{eqn:pairing}
\tUUp \otimes \CV^{\emph{op}} \xrightarrow{\langle \cdot , \cdot \rangle} \BF
\end{equation}
given for all $R \in \CV_{\bn}$ and all $i_1,\dots,i_n \in I$, $d_1,\dots,d_n \in \BZ$ by: 
\begin{equation}
\label{eqn:pairing formula}
\Big \langle \te_{i_1,d_1} \dots \te_{i_n,d_n}, R \Big \rangle = \int_{|z_1| \gg \dots \gg |z_n|} \frac {z_1^{d_1}\dots z_n^{d_n} R(z_1,\dots,z_n)}{\prod_{1\leq a < b \leq n} \zeta_{i_b i_a} \left( \frac {z_b}{z_a} \right)} \prod_{a = 1}^n Dz_a 
\end{equation}
if $\bs_{i_1}+\dots+\bs_{i_n} = \bn$, and 0 otherwise. Implicit in the notation \eqref{eqn:pairing formula} is that the symbol $z_a$ is plugged into one of the variables $z_{i_a \bullet_a}$ of $R$, for all $a \in \{1,\dots,n\}$ \footnote{The choice of $\bullet_a \in \{1,\dots,n_{i_a}\}, \forall a \in \{1,\dots,n\}$ is immaterial due to the symmetry of $R$, as long as we ensure that $\bullet_a \neq \bullet_b$ for all $a \neq b$ such that $i_a = i_b$.}. 

\medskip

\end{proposition}

\begin{remark} While we have presented \eqref{eqn:pairing} as an $\BF$-linear pairing of vector spaces, note that it naturally extends to a bialgebra pairing once one enlarges the algebras involved according to the standard procedure for quantum loop groups (which we recall in Section \ref{sec:double}). The notation $\CV^{\emph{op}}$ in \eqref{eqn:pairing} is meant to underscore this fact.

\end{remark}

\medskip

\begin{proof} To check that \eqref{eqn:pairing} is a well-defined pairing, we need to make sure that any linear relations between the $\te_{i,d}$'s yield linear relations between the right-hand sides of \eqref{eqn:pairing formula}, for any $R \in \CV^{\op}$. Let us first observe that because the defining relations in the algebra $\tUUp$ are quadratic, it suffices to treat the case $n=2$. Indeed, if one inserts in the l.h.s. of \eqref{eqn:pairing formula} an expression $e_{i_1,d_1} \cdots e_{i_{l-1},d_{l-1}} P e_{i_{l+2},d_{l+2} \cdots e_{i_n, d_n}}$ with $P=P(e_{i_l,\bullet},e_{i_{l+1},\bullet})$ a quadratic relation, then for any $R(z_1, \ldots, z_n) \in \CV_n$ we may invoke Fubini's theorem to compute the r.h.s. of \eqref{eqn:pairing formula} by first performing the integration over the variables $z_l,z_{l+1}$.
	
Taking $n=2$, we may rewrite \eqref{eqn:pairing formula} in terms of the generating series $e_i(x)$ and $e_j(y)$ as follows:
$$
\Big \langle \te_i(x) \te_j(y), R \Big \rangle = \int_{|z_1| \gg |z_2|}  \frac {\delta \left(\frac {z_1}x \right)\delta \left(\frac {z_2}y \right) R(z_1,z_2)}{\zeta_{ji} \left( \frac {z_2}{z_1} \right)} Dz_1 Dz_2
$$
where the $\delta$ function is $\delta (x) = \sum_{d \in \BZ} x^d$. Therefore, we also have:
\begin{align}
\Big \langle \te_i(x) \te_j(y) \tzeta_{ji} \left(\frac yx \right) x^{\delta_j^i} , R \Big \rangle = \int_{|z_1| \gg |z_2|}  \frac {\delta \left(\frac {z_1}x \right)\delta \left(\frac {z_2}y \right) \tzeta_{ji} \left(\frac yx \right) x^{\delta_j^i} R(z_1,z_2)}{\zeta_{ji} \left( \frac {z_2}{z_1} \right)} Dz_1 Dz_2 = 
\\[4pt]
\label{eqn:quad pair 1}
= \int_{|z_1| \gg |z_2|}  \frac {\delta \left(\frac {z_1}x \right)\delta \left(\frac {z_2}y \right) \tzeta_{ji} \left(\frac {z_2}{z_1} \right) z_1^{\delta_j^i} R(x,y)}{\zeta_{ji} \left( \frac {z_2}{z_1} \right)} Dz_1 Dz_2   =  R(x,y)(x-y)^{\delta_j^i}
\end{align}
(the equalities above are due to the well-known property: 
\begin{equation}
\label{eqn:delta prop}
\delta \left(\frac zx \right) P(z) = \delta \left(\frac zx \right) P(x)
\end{equation} 
for any Laurent polynomial $P$). Analogously, we have:
\begin{equation}
\label{eqn:quad pair 2}
\Big \langle \te_j(y) \te_i(x)  \tzeta_{ij} \left(\frac xy \right) (-y)^{\delta_j^i}, R \Big \rangle = R(x,y)(x-y)^{\delta_j^i}
\end{equation}
so the two pairings \eqref{eqn:quad pair 1} and \eqref{eqn:quad pair 2} are equal, as we needed to show.
\end{proof}

\medskip

\subsection{} We will now provide a basis of $\tUUp$, following \cite{N} (which was inspired by ideas of \cite{LR, L, R} and \cite{NT}). Consider the set of \textbf{letters}:
$$
i^{(d)} \qquad \forall i \in I, d \in \BZ
$$ 
As usual, a \textbf{word} is any sequence of letters:
$$
\left[ i_1^{(d_1)} \dots i_n^{(d_n)} \right] \qquad \forall i_1,\dots,i_n \in I, d_1, \dots, d_n \in \BZ.
$$
If
$w=\left[ i_1^{(d_1)} \dots i_n^{(d_n)} \right]$ is a word we call $n$ its \textbf{length}, define its \textbf{degree} as:
$$
\deg w = (\bs_{i_1} + \dots + \bs_{i_n}, d_1+\dots+d_n) \in \nn \times \BZ
$$
set
\begin{equation}\label{eqn:defvbar}
\overline{v}=(d_1, \ldots, d_n)
\end{equation}
and put
\begin{equation}
\label{eqn:e w}
\te_w = \te_{i_1,d_1} \dots \te_{i_n,d_n} \in \tUUp
\end{equation}
We write $\mathcal{W}$ for the set of all words. By definition, $\tUUp$ is linearly spanned by the collection of elements $e_w$ for $w \in \mathcal{W}$, but there are linear relations among them. We will now point out a subset of words, such that the corresponding elements yield a linear basis of $\tUUp$. For this, we introduce a total order on the set of words.
We begin by fixing a total order on the set of vertices $I$ of $Q$. \medskip

\begin{definition}[\cite{NT}] 
\label{def:order}

We define a total order on the set of letters as follows:
\begin{equation}
\label{eqn:lex affine}
i^{(d)} < j^{(e)} \quad \text{if} \quad 
\begin{cases} d>e \\ \text{ or } \\ d = e \text{ and } i<j \end{cases}
\end{equation}
We extend this to the total lexicographic order on words by:
$$
\left[ i_1^{(d_1)} \dots i_n^{(d_n)} \right] < \left[ j_1^{(e_1)} \dots j_{m}^{(e_{m})} \right]
$$
if $i_1^{(d_1)} = j_1^{(e_1)}$, \dots, $i_k^{(d_k)} = j_k^{(e_k)}$ and either $i_{k+1}^{(d_{k+1})} < j_{k+1}^{(e_{k+1})}$ or $k = n < m$. \medskip

\end{definition}

\begin{definition}[\cite{N}]
\label{def:non-inc}

A word $w = \left[ i_1^{(d_1)} \dots i_n^{(d_n)} \right]$ is called \textbf{non-increasing} if:
\begin{equation}
\label{eqn:proto}
\begin{cases} d_a < d_b + \sum_{a \leq s < b} \#_{i_si_b} \\ \quad \text{or} \\ d_a = d_b + \sum_{a \leq s < b} \#_{i_si_b} \text{ and } i_a \geq i_b \end{cases}
\end{equation}
for all $1 \leq a < b \leq n$ (see \eqref{eqn:number of edges} for the definition of $\#_{ij}$). \medskip

\end{definition}

\subsection{}

Let $\mathcal{W}_\leq$ denote the set of non-increasing words. The motivation for introducing them is twofold, and is embodied by Lemma~\ref{lem:finite} and Proposition~\ref{prop:proto}.

\medskip

\begin{lemma}[\cite{N}]
\label{lem:finite}

There are finitely many non-increasing words of given degree, which are bounded above by any given word $v$. \medskip

\end{lemma} 

\begin{proof} Let us assume we are counting non-increasing words $[i_1^{(d_1)} \dots i_n^{(d_n)}]$ with $d_1+\dots+d_n = d$ for fixed $n$ and $d$. The fact that such words are bounded above implies that $d_1$ is bounded below. But then the inequality \eqref{eqn:proto} implies that $d_2,\dots,d_n$ are also bounded below. The fact that $d_1+\dots+d_n$ is fixed implies that there can only be finitely many choices for the exponents $d_1,\dots,d_n$. Since there are also finitely many choices for $i_1,\dots,i_n \in I$, this concludes the proof.  
\end{proof} 

\begin{proposition}[\cite{N}]
\label{prop:proto}

The set $\{\te_w\}_{w \in \mathcal{W}_\leq}$ is a linear basis of $\tUUp$. \medskip

\end{proposition}

\begin{proof}[Proof (sketch)] A similar result appears in \cite{N} for a quotient of $\tUUp$. As the proof given in \emph{loc. cit.} only uses relations \eqref{eqn:rel quad}, we may apply it to our situation. More precisely, \loccit shows how to iterate formula \eqref{eqn:rel quad 2} in order to obtain that for any word $v$, the element $\te_v$ belongs to the linear span of elements $\te_w$ with $w \in \mathcal{W}_\leq$ satisfying:
\begin{equation}
\label{eqn:spans}
\begin{split}
w&\geq v,\\
\qquad \min(\bar{v}) - \beta(n) \leq \min(\bar{w}) &\leq \max(\bar{w}) \leq \max(\bar{v}) + \beta(n)
\end{split}
\end{equation}
where 
%$\bar{v} = (d_1,\dots,d_n)$ and
$\beta(n)$ is a universal constant. This implies that:
$$
\tUUp = \text{span} \{\te_w\}_{w \in \mathcal{W}_\leq}
$$
Let us briefly recall how \cite{N} showed that the elements $\{\te_w\}_{w \in \mathcal{W}_\leq}$ are linearly independent, as this will be useful later. Consider any ordered monomial:
\begin{equation}
\label{eqn:monomial}
\mu = z_{i_1\bullet_1}^{-k_1} \dots z_{i_n \bullet_n}^{-k_n}
\end{equation}
where we assume that $\bullet_a \neq \bullet_b$ if $a\neq b$ and $i_a = i_b$. The \textbf{associated word} of $\mu$ is:
\begin{equation}
\label{eqn:leading word}
w_\mu = \left [i_1^{(d_1)} \dots i_n^{(d_n)} \right] 
\end{equation}
where:
\begin{equation}
\label{eqn:k to d}
d_a = k_a + \sum_{t > a} \#_{\overrightarrow{i_ai_t}} - \sum_{s < a} \#_{\overrightarrow{i_si_a}}
\end{equation}
The lexicographically largest of the associated words of various orderings of a given
monomial $\mu$ will be called the \textbf{leading word} of $\mu$ (it is uniquely determined). It was shown in \cite[Lemma 4.8]{N} that the leading word is the only one among all associated words which is non-increasing in the sense of \eqref{eqn:proto}. More generally, the leading word of any non-zero $R \in \CV$, denoted by $\text{lead}(R)$, will be the lexicographically largest of the leading words \eqref{eqn:leading word} for all the monomials which appear in $R$ with non-zero coefficient. Conversely, any $w \in \mathcal{W}_{\leq}$ appears as the leading word:
\begin{equation}
\label{eqn:w as lead}
w = \text{lead}(\text{Sym } \mu)
\end{equation}
of the monomial $\mu$ as in \eqref{eqn:monomial}, chosen such that formula \eqref{eqn:leading word} holds. \medskip

\noindent Analogously to \cite[Formula~(4.18)]{N}, one can show by direct inspection that:
\begin{equation}
\label{eqn:leading word pairing}
\Big \langle e_w, R \Big \rangle \text{ is } \begin{cases} \neq 0 &\text{if }w = \text{lead}(R) \\ = 0 &\text{if }w > \text{lead}(R) \end{cases}
\end{equation}
The formula above immediately shows the linear independence of the elements $e_w$, as $w$ runs over non-increasing words. Indeed, if one were able to write such an element $e_w$ as a linear combination of elements $e_v$ involving only strictly larger non-increasing words $v$ then this would contradict \eqref{eqn:leading word pairing} for $R = \text{Sym }\mu$ with $\mu$ as in \eqref{eqn:w as lead}. 
\end{proof}

\medskip

\begin{corollary}
\label{cor:nondeg}
The pairing $\tUUp \otimes \CV^{\emph{op}} \xrightarrow{\langle \cdot , \cdot \rangle} \BF$ defined in \eqref{eqn:pairing formula} is non-degenerate. \medskip
\end{corollary}

\begin{proof} For non-degeneracy in the second factor, we need to show that any $R \in \CV^{\op}$ which pairs trivially with the whole of $\tUUp$ actually vanishes; this is just the obvious fact that if the power series expansion of the rational function:
$$
\frac {R(z_1,\dots,z_n)}{\prod_{1\leq a < b \leq n} \zeta_{i_bi_a} \left( \frac {z_b}{z_a} \right)}
$$
(in the domain corresponding to an arbitrary order of the variables $z_1,\dots,z_n$ and an arbitrary choice of the indices $i_1,\dots,i_n \in I$) vanishes, then $R = 0$. 

\medskip

\noindent Let us now consider some non-zero element:
$$
\phi = \sum_{w \in \mathcal{W}_\leq} a_w \cdot e_w \in \tUUp
$$
Let $w$ be the smallest word such that $a_w \neq 0$, and choose a monomial $\mu$ whose leading word is $w$ (see \eqref{eqn:monomial}, \eqref{eqn:leading word} and \eqref{eqn:w as lead}). Then by \eqref{eqn:leading word pairing}, $\langle \phi, \text{Sym }\mu\rangle \neq 0$. This gives the non-degeneracy of the pairing in the first factor.
\end{proof}

\medskip

\subsection{}\label{sec:finitedim} For future use, we spell out a ``finite support" variant of Corollary~\ref{cor:nondeg}.  Let $T \subset \mathcal{W}_\leq$ be a finite set of non-increasing words. We set:
$$
\tUU^{+,T} = \bigoplus_{w \in T } \BF \cdot e_w \subset \tUUp
$$
and let $\CV_T \subset \CV$ denote the set of symmetric Laurent polynomials spanned by monomials having the property that their leading word \eqref{eqn:leading word} lies in $T$. Then we claim that the restriction of the pairing \eqref{eqn:pairing} to:
\begin{equation}
\label{eqn:pair rest}
\tUU^{+,T} \otimes \CV^{\op}_T \xrightarrow{\langle \cdot, \cdot \rangle} \BF
\end{equation}
is non-degenerate. Indeed, the two vector spaces above manifestly have the same finite dimension (due to the uniqueness of the leading word \eqref{eqn:leading word} associated to any given monomial), so it suffices to show that \eqref{eqn:pair rest} is non-degenerate in the second argument. This follows from the fact that the leading word $w$ of any $0 \neq R \in \CV_T$ lies in $T$, and \eqref{eqn:leading word pairing} implies that $\langle e_w, R \rangle \neq 0$. 

\bigskip 

\section{The small shuffle algebra and the quantum loop group}
\label{sec:small}

\medskip

\subsection{} We will now define a certain subalgebra of $\CV$, determined by the so-called wheel conditions. These first arose in the context of elliptic quantum groups in \cite{FO}, and the version herein is inspired by the particular wheel conditions of \cite{FHHSY} (which corresponds to the case when $Q$ is the Jordan quiver). \medskip

\begin{definition}[\cite{N}]
\label{def:shuffle}

The \textbf{small shuffle algebra} is the subspace $\CS \subset \CV$ consisting of Laurent polynomials $R(\dots, z_{i1}, \dots, z_{in_i}, \dots)$ that satisfy the  wheel conditions:
\begin{equation}
\label{eqn:wheel}
R \Big|_{z_{ia} = q z_{ic}, z_{jb} = t_e z_{ic}} = 0 
\end{equation}
for any edge $\overline{E} \ni e =  \oij$ and all $a \neq c$, and further $a \neq b \neq c$ if $i = j$. \medskip

\end{definition}

\noindent It is well-known (and straightforward to show) that $\CS$ as defined above is a subalgebra of $\CV$. Because the wheel conditions are vacuous for $R$ of horizontal degree $\{\bs_i\}_{i \in I}$, we conclude that the homomorphism \eqref{eqn:hom big shuffle} actually maps into the small shuffle algebra, i.e.:
\begin{equation}
\label{eqn:big hom shuffle}
\widetilde{\Upsilon}:~\tUUp \longrightarrow \CS 
\end{equation}
The map $\widetilde{\Upsilon}$ was shown to be surjective in \cite{N}. We will obtain in Theorem~\ref{thm:main} below a generators-and-relations description of the algebra $\tUUp / \text{Ker }\widetilde{\Upsilon} \cong \CS$, by describing a set of generators for the kernel of the map $\widetilde{\Upsilon}$. \medskip

\begin{example}
\label{ex:a2}

Let us consider the quiver $Q$ of type $A_2$, consisting of two vertices $\{i,j\}$ with a single edge, say $e\colon i \to j$. Up to gauge transformation, we may specialize the equivariant parameter to $t_{\edge} = q^{\frac 12}$. In this case, it is known that the small shuffle algebra $\CS$ is isomorphic to the positive half of the quantum loop group $U_{q^{-\frac 12}}^+(L\fsl_3)$, which is the quotient of $\tUUp$ by the set of cubic $q$-Serre relations, i.e.:
\begin{equation}
\label{eqn:q-serre}
\CS \simeq \tUUp \Big / \Big( P_{s,t}(x_1,x_2,y) + P_{s,t}(x_2,x_1,y) \Big)
\end{equation}
where for any $s\neq t \in \{i,j\}$, we have set:
\begin{equation*}
\begin{split}
P_{s,t}&(x_1,x_2,y) \\
 \;\; &=e_s(x_1)e_s(x_2)e_t(y) - \left(q^{\frac 12} + q^{-\frac 12} \right) e_s(x_1)e_t(y) e_s(x_2) + e_t(y) e_s(x_1)e_s(x_2). 
\end{split}
\end{equation*}

\end{example}

\medskip

\subsection{} 
\label{sub:cubic}

The wheel conditions \eqref{eqn:wheel} can be interpreted as certain linear conditions on elements $R \in \CV$. As can be expected in light of Corollary~\ref{cor:nondeg}, these linear conditions are given by taking the pairing \eqref{eqn:pairing} with certain elements of $\tUUp$, which we now explicitly describe. We first introduce some more notation. Recall that:
\begin{equation}
\label{eqn:def tzeta}
\tzeta_{ij}(x) = \zeta_{ij}(x) \cdot (1-x)^{\delta_j^i} \in \BF[x^{\pm 1}]
\end{equation}
Given a Laurent polynomial $P(x,y,z)$ and three formal series $e_i(x)$, $e_j(y)$, $e_k(z)$ as in Definition~\ref{def:big quantum}, we will write:
\begin{equation}
\label{eqn:constant term}
\Big[ P(x,y,z) e_i(x) e_j(y) e_k(z) \Big]_{\text{ct}} \in \tUUp
\end{equation}
for the constant term of the expression in square brackets in \eqref{eqn:constant term}. For example, if $P(x) = x^a y^b z^c$ for various integers $a,b,c$, then \eqref{eqn:constant term} equals $e_{i,a}e_{j,b}e_{k,c}$. \medskip

\noindent For any edge $\overline{E} \ni e = \oij$ and any triple of integers $(a,b,c) \in \BZ^3$, we define:
\begin{equation}
\label{eqn:defA}
A^{(e)}_{a,b,c}=\left[ \frac{x_1^a \left( \frac{x_2}{q}\right)^b \left(\frac{y}{t_e}\right)^c }{(1-q)\left(1 - t_{\edge}\right)^{\delta_j^i}\left(1 - \frac {t_{\edge}}q\right)^{\delta_j^i}} \cdot X^{(e)}(x_1,x_2,y) \right]_{\text{ct}} \in \tUU^+_{2\bs_i+\bs_j,a+b+c} 
\end{equation}
where:
\begin{equation}\label{eqn:defX}
\begin{split}
X^{(\edge)}(x_1,x_2,y) = &\frac {\tzeta_{ii}\left(\frac {x_2}{x_1} \right) \tzeta_{ji} \left( \frac {y}{x_1} \right) \tzeta_{ji} \left( \frac {y}{x_2} \right)}{\left(1- \frac {x_2}{x_1q}\right) \left(1-\frac {yq}{x_2t_{\edge}} \right)} \cdot e_i(x_1) e_i(x_2) e_j(y) \\
&+ \frac {\tzeta_{ii}\left(\frac {x_1}{x_2} \right) \tzeta_{ji} \left( \frac y{x_2} \right) \tzeta_{ij} \left( \frac {x_1}y \right) \left(- \frac {x_2 t_{\edge}}{y} \right) \left(- \frac y{x_1} \right)^{\delta_j^i} }{\left(1- \frac {yq}{x_2t_{\edge}}\right) \left(1-\frac {x_1 t_{\edge}}{y} \right)} \cdot e_i(x_2) e_j(y) e_i(x_1) \\
&+ \frac {\tzeta_{ii}\left(\frac {x_2}{x_1} \right) \tzeta_{ij} \left( \frac {x_1}y \right) \tzeta_{ij} \left( \frac {x_2}y \right) \left( \frac {x_2t_{\edge}}{yq} \right)\left( \frac {y^2}{x_1x_2} \right)^{\delta_j^i}}{\left(1- \frac {x_2}{x_1 q}\right) \left(1-\frac {x_1 t_{\edge}}{y} \right)} \cdot e_j(y) e_i(x_1) e_i(x_2)  
\end{split}
\end{equation}
The linear factors in $x_1,x_2,y$ in the denominators above are all canceled by some factors in the numerator, hence the expression in the square brackets of \eqref{eqn:defA} is a Laurent polynomial in $x_1,x_2,y$ times the product of the series $e_i(x_1)$, $e_i(x_2)$, $e_j(y)$.

\medskip

\begin{proposition}
\label{prop:wheel pairing}
The elements $A^{(e)}_{a,b,c}$ only depend on $d=a+b+c$ (and will henceforth simply be denoted $A^{(e)}_{d}$). Setting:  
$$
A^{(\edge)}(x) = \sum_{d \in \BZ} A^{(\edge)}_dx^{-d} 
$$
we have:
\begin{equation}
\label{eqn:wheel pair}
\Big \langle A^{(\edge)}(x), R \Big \rangle = R \Big|_{z_{i1} = x, z_{j1} = t_{\edge} x, z_{i2} = qx} \qquad \forall R \in \CV^{\emph{op}}_{2\bs_i+\bs_j} 
\end{equation}
for any edge $\overline{E} \ni e = \oij$. \medskip

\end{proposition} 

\begin{proof} We need to show that for any integers $a+b+c = d$, we have:
\begin{equation}
\label{eqn:wheel coeff}
\Big \langle A^{(\edge)}_d, R \Big \rangle  = \text{coefficient of }x^{-d} \text{ in } R \Big|_{z_{i1} = x, z_{j1} = t_{\edge} x, z_{i2} = qx}
\end{equation}
for all $R \in \CV^{\op}$. This will also imply that $A^{(\edge)}_{a,b,c}$ only depends on $d$, because of the non-degeneracy of the pairing (Corollary \ref{cor:nondeg}). We have:
$$
A^{(e)}_d = C_1+C_2+C_3
$$
where $C_1,C_2, C_3$ correspond to the three terms defining $X^{(e)}_{a,b,c}$. Abbreviating:
$$\mathbf{D}= \frac {(z_{i1})^a \left(\frac {z_{i2}}q \right)^b \left( \frac {z_{j1}}{t_{\edge}} \right)^c }{(1-q)\left(1 - t_{\edge}\right)^{\delta_j^i}\left(1 - \frac {t_{\edge}}q\right)^{\delta_j^i}}Dz_{i1} Dz_{i2} Dz_{j1}$$
we have, using \eqref{eqn:pairing formula}:
\begin{align*}
&\Big \langle C_1, R \Big \rangle = \int_{|z_{i1}| \gg |z_{i2}| \gg |z_{j1}|}  \frac {R(z_{i1},z_{i2},z_{j1})  \left(1- \frac {z_{i2}}{z_{i1}} \right)\left(1- \frac {z_{j1}}{z_{i1}} \right)^{\delta_j^i} \left(1- \frac {z_{j1}}{z_{i2}} \right)^{\delta_j^i} \mathbf{D}}{\left(1- \frac {z_{i2}}{z_{i1}q}\right) \left(1-\frac {z_{j1}q}{z_{i2}t_{\edge}} \right)}  \\
&\Big \langle C_2, R \Big \rangle = \int_{|z_{i2}| \gg |z_{j1}| \gg |z_{i1}|} \frac {R(z_{i1},z_{i2},z_{j1})\left(\frac {z_{i1}t_{\edge}}{z_{j1}} \right)\left(1- \frac {z_{i2}}{z_{i1}} \right)\left(1- \frac {z_{j1}}{z_{i1}} \right)^{\delta_j^i}\left(1- \frac {z_{j1}}{z_{i2}} \right)^{\delta_j^i}\mathbf{D}}{\left(1- \frac {z_{j1}q}{z_{i2}t_{\edge}}\right) \left(1-\frac {z_{i1}t_{\edge}}{z_{j1}} \right)}  \\
&\Big \langle C_3, R \Big \rangle = \int_{|z_{j1}| \gg |z_{i1}| \gg |z_{i2}|} \frac {R(z_{i1},z_{i2},z_{j1}) \left(\frac {z_{i2}}{z_{i1}q}\right) \left(\frac {z_{i1}t_{\edge}}{z_{j1}} \right)\left(1- \frac {z_{i2}}{z_{i1}} \right)\left(1- \frac {z_{j1}}{z_{i1}} \right)^{\delta_j^i}\left(1- \frac {z_{j1}}{z_{i2}} \right)^{\delta_j^i}\mathbf{D}}{\left(1- \frac {z_{i2}}{z_{i1}q}\right) \left(1-\frac {z_{i1}t_{\edge}}{z_{j1}} \right)} 
\end{align*}
(we replace $z_{j1}$ by $z_{i3}$ in the formulas above if $i = j$; note that in the middle formula for $\langle C_2, R \rangle$, we took the liberty of replacing $z_{i1} \leftrightarrow z_{i2}$ in the arguments of the symmetric Laurent polynomial $R$). Then \eqref{eqn:wheel coeff} is an immediate consequence of the following identity of formal series:
\begin{align}
\delta \left(\frac {z_{i2}}{z_{i1}q} \right)\delta \left(\frac {z_{j1}}{z_{i1}t_{\edge}} \right) &= ev_{ |z_{i1}| \gg |z_{i2}| \gg |z_{j1}| } \left[\frac 1{\left(1- \frac {z_{i2}}{z_{i1}q}\right) \left(1-\frac {z_{j1}q}{z_{i2}t_{\edge}} \right)} \right] \\
&+  ev_{ |z_{i2}| \gg |z_{j1}| \gg |z_{i1}| } \left[\frac {\left(\frac {z_{i1}t_{\edge}}{z_{j1}} \right)}{\left(1- \frac {z_{j1}q}{z_{i2}t_{\edge}}\right) \left(1-\frac {z_{i1}t_{\edge}}{z_{j1}} \right)}\right]  \\
&+  ev_{ |z_{j1}| \gg |z_{i1}| \gg |z_{i2}| } \left[\frac {\left(\frac {z_{i2}}{z_{i1}q}\right) \left(\frac {z_{i1}t_{\edge}}{z_{j1}} \right)}{\left(1- \frac {z_{i2}}{z_{i1}q}\right) \left(1-\frac {z_{i1}t_{\edge}}{z_{j1}} \right)}\right] \label{eqn:delta functions}
\end{align}
where $ev_{|x| \gg |y| \gg |z|}$ is the operator of taking the Laurent expansion of a rational function in the asymptotic region $ |x| \gg |y| \gg |z|$. The formula above is a straightforward computation involving formal power series, which we leave as an exercise to the interested reader. 
\end{proof}

\medskip

\subsection{}\label{sec:3.7}

By Proposition~\ref{prop:wheel pairing}, different choices of integers $a,b,c$ yield elements $A^{(\edge)}_d$ which are equal modulo the ideal generated by relations \eqref{eqn:rel quad}. In order to write $A^{(\edge)}_d$  ``canonically", one can expand these elements in the basis $\{e_w\}_{w \in \mathcal{W}_\leq}$ of $\tUUp$ (recall Proposition \ref{prop:proto}):
\begin{equation}
\label{eqn:expand non-inc}
A^{(\edge)}_d =\sum_{w = [ i_1^{(\alpha_1)} i_2^{(\alpha_2)} i_3^{(\alpha_3)} ] \in \mathcal{W}_\leq} c^{(e)}_{d,w}\cdot e_{i_1,\alpha_1} e_{i_2,\alpha_2} e_{i_3,\alpha_3}
\end{equation}
The coefficients $c_{d,w}^{(e)}$ occurring in the expression above are not nice (even in relatively simple cases, such as a vertex with a single loop, or two distinct vertices with two edges between them), but it is easy to see that the transformation: 
$$
d\mapsto d+3\ , \qquad w = [i_1^{(\alpha_1)}i_2^{(\alpha_2)}i_3^{(\alpha_3)}] \mapsto [i_1^{(\alpha_1+1)}i_2^{(\alpha_2+1)}i_3^{(\alpha_3+1)}].
$$
rescales $c_{d,w}^{(e)}$ by a monomial in $q,t_e$. This allows us to conclude that there exists a large enough natural number $N$, which only depends on $Q$, such that:
\begin{equation}
\label{eqn:large enough}
|\alpha_t - \alpha_s| \leq N, \ \forall s,t \in  \{1,2,3\}
\end{equation}
for all words $[i_1^{(\alpha_1)}i_2^{(\alpha_2)}i_3^{(\alpha_3)}]$ which appear with non-zero coefficient in \eqref{eqn:expand non-inc}, for all edges $\edge \in E$ and all $d \in \BZ$. We will use this fact in the proof of Theorem~\ref{thm:main}. \medskip

\subsection{} Our main result is that the two-sided ideal $J \subset \tUUp$ generated by $A^{(\edge)}_d$, for all $\edge \in \overline{E}$ and $d \in \BZ$, coincides with the kernel of the map \eqref{eqn:big hom shuffle}.  
To this end, we define the (positive half of the) \textbf{generic quantum loop group} as:
\begin{equation}
\label{eqn:small quantum}
\UUp = \tUUp \Big / J.
\end{equation}
Note that taking the quotient by $J$ is equivalent to imposing $X^{(\edge)}(x_1,x_2,y) = 0$ for all edges $\edge$, which is precisely the content of \eqref{eqn:serre a}. Our main result, which is simply a restatement of Theorem \ref{thm:intro}, is the following. \medskip

\begin{theorem} 
\label{thm:main}
The assignment $\te_{i,d} \mapsto z_{i1}^d$ for $i \in I, d \in \BZ$ induces an algebra isomorphism:
$$
\Upsilon : \UUp \xrightarrow{\sim} \CS.
$$

\medskip

\end{theorem}

\begin{example}
\label{ex:a1}

Let us work out the cubic relation \eqref{eqn:serre a} when $Q$ is the $A_2$ quiver of Example~\ref{ex:a2}, with the goal of recovering the $q$-Serre relation \eqref{eqn:q-serre}. Recall that we set $t_{\edge} = q^{\frac 12}$, so relation \eqref{eqn:serre a} states that:
\begin{equation}\label{eqn:a1}
\begin{split}
\left(x_1 - y q^{\frac 12}\right) e_i(x_1) e_i(x_2) e_j(y) + &\left(x_2q^{\frac 12} - x_1 q^{-\frac 12} \right) e_i(x_2) e_j(y) e_i(x_1) \\
 &+ \left(yq^{-\frac 12} - x_2 \right) e_j(y) e_i(x_1) e_i(x_2) = 0
\end{split}
\end{equation}
The quadratic relations \eqref{eqn:rel quad} hold in both $\UUp$ and $U_{q^{-\frac 12}}^+(L\fsl_3)$, and they read:
\begin{align}
e_s(x)e_s(y) \left(x q - y \right) &= e_s(y)e_s(x)  \left(x - y q\right) \label{eqn:a2} \\
e_s(x)e_t(y) \left(x - y q^{\frac 12} \right) &= e_t(y)e_s(x) \left(x q^{\frac 12} - y\right) \label{eqn:a3}
\end{align}
for $s\neq t \in \{i,j\}$. Let us show how to obtain the $q$-Serre relation (for $s=i$, $t=j$) from \eqref{eqn:a1}, \eqref{eqn:a2}, \eqref{eqn:a3}. As a consequence of \eqref{eqn:a3}, we have:
\begin{align*}
\left(x_2 q^{\frac 12} - y q \right) e_i(x_2)e_j(y)e_i(x_1) &= \left(x_2q - yq^{\frac 12} \right) e_j(y) e_i(x_2) e_i(x_1) \\
\left(y q^{-1} - x_1 q^{-\frac 12} \right) e_i(x_2)e_j(y)e_i(x_1) &= \left(yq^{-\frac 12} - x_1q^{-1} \right) e_i(x_2) e_i(x_1)  e_j(y)
\end{align*}
Adding the two relations together, and then subtracting \eqref{eqn:a1} from the result yields:
\begin{equation*}
\begin{split}
y \left(q^{-1} - q \right)& e_i(x_2)e_j(y)e_i(x_1) \\ 
=& \left(x_2q - yq^{\frac 12} \right) e_j(y) e_i(x_2) e_i(x_1) + \left(yq^{-\frac 12} - x_1q^{-1} \right) e_i(x_2) e_i(x_1)  e_j(y)\\
&+ \left(x_1 - y q^{\frac 12} \right) e_i(x_1) e_i(x_2) e_j(y) + \left(yq^{-\frac 12} - x_2 \right) e_j(y) e_i(x_1) e_i(x_2)
\end{split}
\end{equation*}
Symmetrizing the relation above with respect to $x_1 \leftrightarrow x_2$ yields:
\begin{equation*}
\begin{split}
y \left(q^{-1} - q \right) &\Big[ e_i(x_2)e_j(y)e_i(x_1) + e_i(x_1)e_j(y)e_i(x_2) \Big] = \\
&= \left[x_2q - yq^{\frac 12} + yq^{-\frac 12} -  x_1 \right] e_j(y) e_i(x_2) e_i(x_1) \\
&+ \left[ yq^{-\frac 12} - x_1q^{-1} + x_2 - yq^{\frac 12} \right] e_i(x_2) e_i(x_1)  e_j(y) \\
&+ \left[x_1 - y q^{\frac 12} + yq^{-\frac 12} - x_2q^{-1} \right] e_i(x_1) e_i(x_2) e_j(y) \\
&+ \left[yq^{-\frac 12}  - x_2 + x_1q - y q^{\frac 12} \right] e_j(y) e_i(x_1) e_i(x_2)
\end{split}
\end{equation*}
By applying \eqref{eqn:a2}, the right-hand side of the expression above is equal to:
\begin{equation*}
\begin{split}
y \left(q^{-\frac 12} - q^{\frac 12} \right)\Big[ &e_j(y) e_i(x_2)e_i(x_1) + e_i(x_2)e_i(x_1)e_j(y) \\
&+ e_i(x_1)e_i(x_2)e_j(y) + e_j(y) e_i(x_1)e_i(x_2) \Big]
\end{split}
\end{equation*}
Dividing by $y(q^{-\frac 12}-q^{\frac 12})$, we obtain precisely the $q$-Serre relation: 
\begin{equation}
\label{eqn:cubic explicit}
P_{s,t}(x_1,x_2,y) + P_{s,t}(x_2,x_1,y) = 0
\end{equation}
of \eqref{eqn:q-serre}. We have just showed that the usual $q$-Serre relations hold in the generic quantum loop group. The distinction between the cubic relations \eqref{eqn:cubic explicit} and the cubic relations \eqref{eqn:serre a} is explained by the fact that the two sets of relations can be obtained from each other by adding appropriate multiples of the quadratic relations. \medskip
\end{example}

\begin{proof}[Proof of Theorem~\ref{thm:main}] Let us recall the non-degenerate pairing:
\begin{equation}
\label{eqn:pairing shuffle}
\CS \otimes \CS^{\op} \xrightarrow{\langle \cdot , \cdot \rangle'} \BF
\end{equation}
defined in Proposition 3.3 of \cite{N}. Comparing our formula \eqref{eqn:pairing formula} with formula (3.30) of \loccitt, we see that the two pairings are compatible, in the sense that:
\begin{equation}
\label{eqn:compatible pairings}
\Big \langle f,\iota(g) \Big \rangle = \Big \langle \widetilde{\Upsilon}(f), g\Big \rangle' \in \mathbb{F}, \qquad \forall\; f \in \tUUp, g \in \CS^{\op},
\end{equation}
where $\iota : \CS^{\op} \to \CV^{\op}$ is the tautological inclusion.
To show that the homomorphism $\widetilde{\Upsilon}$ of \eqref{eqn:big hom shuffle} descends to a homomorphism:
\begin{equation}
\label{eqn:small hom shuffle}
\Upsilon:\UUp \rightarrow \CS,
\end{equation}
one needs to prove that for any $e \in \overline{E}$ and $d \in \BZ$, we have:
\begin{equation}
\label{eqn:elements}
\widetilde{\Upsilon}(A^{(e)}_d) = 0.
\end{equation}
By the non-degeneracy of the pairing \eqref{eqn:pairing shuffle}, it suffices to show that $\widetilde{\Upsilon}(A^{(e)}_d)$ pair trivially with anything in $\CS^{\op}$. Using \eqref{eqn:compatible pairings}, this is equivalent to showing that $A^{(e)}_d$ pair trivially with anything in $\CS^{\op}$ under the pairing \eqref{eqn:pairing}, which follows from \eqref{eqn:wheel pair}. \medskip

\noindent Since the homomorphism $\widetilde{\Upsilon}$ of \eqref{eqn:big hom shuffle} is surjective, so is $\Upsilon$. It thus remains to show that $\Upsilon$ is injective, and because the pairing \eqref{eqn:pairing shuffle} is non-degenerate, the task boils down to showing that an element:
\begin{equation}
\label{eqn:relation}
\psi = \sum_{w \in\mathcal{W}_\leq} c_w \cdot \te_w \in \tUUp
\end{equation}
which pairs trivially with $\CS^{\op}$ lies in $J$. Clearly, it suffices to do so for a homogeneous $\psi$. Let us fix a degree $(\bn,d) \in \nn \times \BZ$, fix some $\psi \in \tUUp$ of degree $(\bn,d)$ as above and define, for a pair of positive integers $M,m$:
\begin{multline}
\label{eqn:def t}
T_{M,m} = \Big\{ w =\left[i_1^{(d_1)} \dots i_n^{(d_n)} \right]  \in \mathcal{W}_\leq,\;  \deg(w)= (\bn,d), \\ \sum_{s \in A} d_s \geq -  M |A| + m |A|^2 - \sum_{A \ni s > t \notin A} \#_{i_si_t}, \quad \forall A \subseteq \{1,\dots,n\} \Big\}.
\end{multline}
For notational simplicity, we will write $T=T_{M,m}$. Observe that $T$ is finite, and that if $M$ is picked large enough, the set $T$ will contain all the words which appear with non-zero coefficient in \eqref{eqn:relation}. As we are free to pick $M$ and $m$ arbitrarily large, it might seem strange that we bother to subtract the quantity $\sum_{A \ni s > t \notin A} \#_{i_si_t}$ from the right-hand side of the inequality in \eqref{eqn:def t}. This is justified by the following claim, whose straightforward proof we leave to the reader:

\begin{claim}
	 The set of inequalities 
$$\sum_{s \in A} d_s \geq -  M |A| + m |A|^2 - \sum_{A \ni s > t \notin A} \#_{i_si_t}, \quad \forall A \subseteq \{1,\dots,n\}$$
holds for the leading word $\left[i_1^{(d_1)} \dots i_n^{(d_n)} \right]$ of a monomial $\mu$ if and only if it holds for the associated word \eqref{eqn:leading word} of any other ordering of the variables of $\mu$.
\end{claim}

With this in mind, we have reduced the proof of Theorem~\ref{thm:main} to the following lemma: \medskip

\begin{lemma}
\label{claim:1}

Let $\CS_{T}^{\emph{op}} = \CS^{\emph{op}} \cap \CV_{T}^{\emph{op}}$. For $M \gg m \gg 1$, the pairing:
\begin{equation}
\label{eqn:pairing in claim}
\left( \tUU^{+,T} \Big/ J \cap \tUU^{+,T} \right) \otimes \CS^{\emph{op}}_{T} \xrightarrow{\langle \cdot, \cdot \rangle_T'} \BF
\end{equation}
induced by \eqref{eqn:pair rest} is non-degenerate. \medskip

\end{lemma}

\begin{proof}[Proof of Lemma~\ref{claim:1}] Let us denote by $\langle\cdot,\cdot\rangle_T$ the restriction of $\langle\cdot,\cdot\rangle$ to $\tUU^{+,T} \otimes \CV^{\op}_T$, as in \eqref{eqn:pair rest}. By construction, the pairings $\langle \cdot, \cdot \rangle'_T$ and $\langle\cdot,\cdot \rangle_T$ are compatible in the sense that:
$$
\Big \langle f,\iota_T(g) \Big \rangle_T = \Big \langle \widetilde{\Upsilon}(f), g\Big \rangle_T' \in \mathbb{F}, \qquad \forall\; f \in \tUU^{+,T}, g \in \CS^{\op}_T,
$$
where $\iota : \CS^{\op}_T \to \CV^{\op}_T$ is the tautological inclusion; observe that $\widetilde{\Upsilon}( \tUU^{+,T}) \subseteq \CS_{T}$. By Section~\ref{sec:finitedim}, $\langle \cdot , \cdot \rangle_T$ is a non-degenerate pairing of \textit{finite-dimensional} vector spaces. To show that $\langle \cdot , \cdot \rangle_T'$ is also non-degenerate, it hence suffices to check that:
\begin{equation}
\label{eqn:suffices 1}
\CS_{T}^{\op} = (J \cap \tUU^{+,T})^\perp 
\end{equation}
Note that the inclusion $\CS_T^{\op} \subseteq (J \cap \tUU^{+,T})^\perp$ is obvious, so we will prove the opposite:
\begin{equation}
\label{eqn:suffices 2}
\CS_T^{\op} \supseteq (J \cap \tUU^{+,T})^\perp
\end{equation}
Equivalently, this requires us to show that:
\begin{equation}
\label{eqn:suffices 3}
\forall R \in \CV_{T}^{\op} \backslash \CS_{T}^{\op}, \qquad \exists \psi \in J \cap \tUU^{+,T} \quad \text{such that} \quad \langle \psi, R \rangle \neq 0
\end{equation}
Let us fix $R$ as above. As it does not satisfy one of the wheel conditions, we have:
\begin{equation}
\label{eqn:edge non-zero}
R \Big|_{z_{ia} =  qx, z_{jb} = t_{\edge} x, z_{ic} = x} \neq 0
\end{equation}
for some edge $e = \oij$. By \eqref{eqn:pairing formula} and \eqref{eqn:wheel pair}, this entails:
\begin{equation}
\label{eqn:not zero}
\left \langle e_{i_1,d_1} \dots e_{i_{n-3},d_{n-3}} A^{(\edge)}_{d}, R \right \rangle \neq 0
\end{equation}
for certain $i_1,\dots,i_{n-3} \in I$ and $d_1,\dots,d_{n-3},d \in \BZ$ (to see this claim, one must expand the integral \eqref{eqn:pairing formula} in the domain where the variables $z_{ia}, z_{jb}, z_{ic}$ of \eqref{eqn:edge non-zero} have smaller absolute value than all other variables of $R$). Assume that the word:
\begin{equation}
\label{eqn:word not zero}
\left[ i_1^{(d_1)} \dots i_{n-3}^{(d_{n-3})} \right]
\end{equation}
is maximal such that \eqref{eqn:not zero} holds. In particular, this implies that the word \eqref{eqn:word not zero} must be non-increasing (as we have seen in the proof of Proposition~\ref{prop:proto}, any $e_v$ can be written as a linear combination of $e_w$ as $w > v$ runs over non-increasing words). To establish the required \eqref{eqn:suffices 3}, it therefore suffices to show that:
\begin{equation}
\label{eqn:suffices 4}
\psi := e_{i_1,d_1} \dots e_{i_{n-3},d_{n-3}} A^{(\edge)}_{d} \in \tUU^{+,T}
\end{equation}
To this end, let us label the variables of $R$ other than $z_{ia},z_{jb},z_{ic}$ by $z_1,\dots,z_{n-3}$ (of colors $i_1,\dots,i_{n-3}$, respectively). Using formula \eqref{eqn:pairing formula}, relation \eqref{eqn:not zero} states:
\begin{multline}
\int_{|z_1| \gg \dots \gg |z_{n-3}| \gg |x|}  Dx \prod_{a = 1}^{n-3} Dz_a \\ \frac {z_1^{d_1}\dots z_{n-3}^{d_{n-3}} x^d R(z_1,\dots,z_{n-3},qx,t_ex,x)}{\prod_{1 \leq a < b \leq n-3} \zeta_{i_b i_a} \left( \frac {z_b}{z_a} \right) \prod_{a=1}^{n-3} \left[ \zeta_{i i_a} \left( \frac {qx}{z_a} \right) \zeta_{j i_a} \left( \frac {t_ex}{z_a} \right) \zeta_{i i_a} \left( \frac {x}{z_a} \right)\right]} \neq 0 \label{eqn:non vanishing}
\end{multline}
As: 
$$
\zeta_{ij}(u) \in u^{-\#_{\oji}} \BF[[u]]^\times, \quad \forall i,j \in I
$$ 
the non-vanishing \eqref{eqn:non vanishing} and the maximality of the word \eqref{eqn:word not zero} imply that:
$$
\left[z_1^{k_1} \dots z_{n-3}^{k_{n-3}} x^{d + \sum_{a=1}^{n-3} \left( 2 \#_{\overrightarrow{i_ai}} + \#_{\overrightarrow{i_aj}}\right)} R(z_1,\dots,z_{n-3},qx,t_ex,x)\right]_{\text{ct}} \neq 0
$$
where we write for all $a \in \{1,\dots,n-3\}$:
$$
k_a = d_a + \sum_{s < a} \#_{\overrightarrow{i_si_a}} - \sum_{t > a} \#_{\overrightarrow{i_ai_t}} - 2 \#_{\overrightarrow{i_ai}} -  \#_{\overrightarrow{i_aj}}
$$
(compare with \eqref{eqn:k to d}). Thus, we conclude that the Laurent polynomial $R \in \CV_T^{\op}$ must include with non-zero coefficient a monomial of the form:
$$
z_1^{-k_1} \dots z_{n-3}^{-k_{n-3}} z_{ia}^{-k_{n-2}} z_{jb}^{-k_{n-1}} z_{jc}^{-k_n}
$$
where:
\begin{align*}
&k_{n-2} = d_{n-2} + \sum_{a = 1}^{n-3} \#_{\overrightarrow{i_ai}} -\#_{\overrightarrow{ij}} -\#_{\overrightarrow{ii}} \\
&k_{n-1} = d_{n-1} + \sum_{a = 1}^{n-3} \#_{\overrightarrow{i_aj}} + \#_{\overrightarrow{ij}} -\#_{\overrightarrow{ji}} \\
&k_n = d_n + \sum_{a = 1}^{n-3} \#_{\overrightarrow{i_ai}} + \#_{\overrightarrow{ji}} + \#_{\overrightarrow{ii}} 
\end{align*}
for some integers $d_{n-2}+d_{n-1}+d_n = d$. We will henceforth write $i_{n-2} = i$, $i_{n-1} = j$ and $i_n = i$. As explained in the paragraph following \eqref{eqn:def t}, although the word:
$$
\left[i_1^{(d_1)} \dots i_{n-3}^{(d_{n-3})} i^{(d_{n-2})} j^{(d_{n-1})} i^{(d_n)} \right]
$$
might fail to be non-increasing (our assumption only states that its prefix obtained by removing the last three letters is non-increasing), its exponents still satisfy the inequality in \eqref{eqn:def t}. In particular, we have:
\begin{equation}
\label{eqn:holds}
\sum_{s \in A} d_s \geq -  M |A| + m |A|^2 - \sum_{A \ni s > t \notin A} \#_{i_si_t}  \end{equation}
where $A = B$ or $A = B \sqcup \{n-2,n-1,n\}$, for arbitrary $B \subseteq \{1,\dots,n-3\}$. However, as explained in Subsection~\ref{sec:3.7}, the element $\psi$ of \eqref{eqn:suffices 4} is a linear combination of products of the form:
\begin{equation}
\label{eqn:product}
e_{i_1,d_1} \dots e_{i_{n-3},d_{n-3}} e_{j_1,\delta_1} e_{j_2,\delta_2} e_{j_3,\delta_3} 
\end{equation}
where $\{j_1,j_2,j_3\} = \{i,i,j\}$ and the numbers $\delta_1,\delta_2,\delta_3$ all satisfy:
\begin{equation}
\label{eqn:c1}
\delta_1+\delta_2+\delta_3 = d \qquad \text{and} \qquad \left|\delta_u - \frac d3 \right| \leq c_1, \forall u \in \{1,2,3\}
\end{equation}
for some global constant $c_1$. For every term of the form \eqref{eqn:product} that appears in $\psi$, let us consider the largest index $x \in \{0,\dots,n-3\}$ such that:
\begin{equation}
\label{eqn:final ineq}
d_x < d_{x+1} - c_2
\end{equation}
where the global constant $c_2$ is chosen much larger than the number $\beta(n)$ that appears in \eqref{eqn:spans}. As explained in Proposition~\ref{prop:proto}, we may write:
$$
e_{i_{x+1},d_{x+1}} \dots e_{i_{n-3},d_{n-3}} e_{j_1,\delta_1} e_{j_2,\delta_2} e_{j_3,\delta_3}  \in \sum_{\left[j_{x+1}^{(l_{x+1})} \dots j_n^{(l_n)} \right] \in \CW_{\leq}} \BF \cdot e_{j_{x+1},l_{x+1}} \dots e_{j_n,l_n} 
$$
in such a way that the exponents $l_{x+1},\dots,l_n$ that appear with non-zero coefficient in the right-hand side satisfy the following properties:

\medskip

\begin{itemize}[leftmargin=*]

\item all the numbers $l_{x+1},\dots,l_n$ are within a global constant $c_3$ away from their average, which will be denoted by $o$,

\medskip

\item the large difference between $d_x$ and $d_{x+1}$ ensures that all concatenated words
\begin{equation}
\label{eqn:concatenated words}
\left[ i_1^{(d_1)} \dots i_x^{(d_x)} j_{x+1}^{(l_{x+1})} \dots j_n^{(l_n)} \right]
\end{equation}
which arise in the procedure above are non-increasing (recall that the word \eqref{eqn:word not zero} was non-increasing to begin with, and thus so are all of its prefixes).

\end{itemize}

\medskip

\noindent To prove \eqref{eqn:suffices 4}, it therefore remains to show that the concatenated words that appear in \eqref{eqn:concatenated words} are in $T$. Property \eqref{eqn:holds} implies that
\begin{align}
&\sum_{s \in B} d_s \geq -M|A| + m|A|^2 - 2n^2|E| \label{eqn:last 1} \\
&\sum_{s \in B} d_s + \sum_{s=x+1}^n l_s \geq -M|A| + m|A|^2 - 2n^2|E| \label{eqn:last 2}
\end{align}
for $A = B$ or $A = B \sqcup \{x+1,\dots,n\}$ respectively, where $B \subseteq \{1,\dots,x\}$ is arbitrary. Assume for the purpose of contradiction that the defining property of $T$ is violated for
$$
A = B \sqcup C
$$
with $B \subseteq \{1,\dots,x\}$ and $C$ a proper subset of $\{x+1,\dots,n\}$, i.e.
\begin{equation}
\label{eqn:last 3}
\sum_{s \in B} d_s + \sum_{s \in C} l_s < -M|A| + m|A|^2 
\end{equation}
We claim that \eqref{eqn:last 1}--\eqref{eqn:last 2} and \eqref{eqn:last 3} are incompatible (for $m$ chosen large enough compared to the constant $c_3$ mentioned in the first bullet above). Indeed, the properties listed in the first bullet above allow us to obtain the following inequalities from \eqref{eqn:last 1}, \eqref{eqn:last 2}, \eqref{eqn:last 3}:
\begin{align}
&\sum_{s \in B} d_s \geq -M|B| + m|B|^2 - 2n^2|E| \label{eqn:last 4} \\
&\sum_{s \in B} d_s + y(o-c_3) < -M(|B|+y) + m(|B|+y)^2  \label{eqn:last 5} \\
&\sum_{s \in B} d_s + (n-x)(o+c_3) \geq -M(|B|+n-x) + m(|B|+n-x)^2 -  2n^2|E| \label{eqn:last 6}
\end{align}
where $y = |C|$ lies in $\{1,\dots,n-x-1\}$. Subtracting \eqref{eqn:last 4} from \eqref{eqn:last 5} yields
$$
o-c_3 < -M + m(y+2|B|) + 2n^2|E|
$$
and subtracting \eqref{eqn:last 5} from \eqref{eqn:last 6} yields 
$$
o+c_3 \geq -M + m(n-x+y + 2|B|) - 2n^2|E|
$$
The two inequalities above are incompatible if $m$ is chosen large enough compared to $c_3$ and $2n^2|E|$, thus yielding the desired contradiction. We have thus proved Lemma~\ref{claim:1}, hence also Theorem~\ref{thm:main}.
\end{proof}
\end{proof}

\bigskip

\section{The Drinfeld double of $\UUp$}\label{sec:double}

\medskip

\subsection{} We will now recall the standard procedure of upgrading $\UUp \cong \CS$ to a Hopf algebra isomorphism, by extending and doubling the two algebras involved (we will follow the conventions of \cite{N}). Consider the opposite algebra:
$$
\UUm = \UU^{+,\op} \cong \CS^{\op}
$$
whose generators will be denoted by $\{f_{i,d}\}_{i\in I, d \in \BZ}$; we set $f_i(z) = \sum_{d \in \BZ} f_{i,d}z^{-d}$. \medskip

\begin{definition}
\label{def:extended}
The extended algebras $\UU^{\geq}$ and $\UU^{\leq}$ are defined as the semi-direct products:
\begin{align*}
&\UU^{\geq} = \UUp \bigotimes_{\BF} \BF [h^+_{i,d}]_{i \in I, d \geq 0} \\
&\UU^{\leq} = \UUm \bigotimes_{\BF} \BF [h^-_{i,d}]_{i \in I, d \geq 0}
\end{align*}
where the multiplication is governed by the following relations for all $i,j \in I$:
\begin{align}
&e_i(z) h^+_j(w) = h^+_j(w) e_i(z) \frac {\zeta_{ij} \left( \frac z{w} \right)}{\zeta_{ji} \left( \frac {w}z \right)} \label{eqn:comm ext 1} \\
&f_i(z) h^-_j(w) = h^-_j(w) f_i(z) \frac {\zeta_{ji} \left( \frac {w}z \right)}{\zeta_{ij} \left( \frac z{w} \right)}. \label{eqn:comm ext 2}
\end{align}
Here $h_j^\pm(w) = \sum_{d = 0}^{\infty} h_{j,d}^\pm w^{\mp d}$. The rational functions in the right-hand sides of \eqref{eqn:comm ext 1} and \eqref{eqn:comm ext 2} are expanded in the same asymptotic direction of $w$ as $h_j^\pm(w)$. \medskip

\end{definition}

\noindent The algebras $\UU^{\geq}$ and $\UU^{\leq}$ are actually bialgebras, with respect to the coproduct:
\begin{equation}
\label{eqn:cop 1}
\Delta \left(h^{\pm}_i(z) \right) = h^{\pm}_i(z) \otimes h^{\pm}_i(z), 
\end{equation}
\begin{equation}
\label{eqn:cop 2}
\Delta \left(e_i(z) \right) = e_i(z) \otimes 1 + h^+_i(z) \otimes e_i(z) ,
\end{equation}
\begin{equation}
\label{eqn:cop 3}
\Delta \left(f_i(z) \right) = f_i(z) \otimes h^-_i(z) + 1 \otimes f_i(z)
\end{equation}
(strictly speaking, the coproduct above consists of infinite sums, meaning that $\UU^{\geq}$ and $\UU^{\leq}$ are topological bialgebras; we will not dwell upon this fact, as it is quite routine in the theory of quantum affinizations). The final step in the construction of the Drinfeld double is to note that there is a (unique) bialgebra pairing:
\begin{equation}
\label{eqn:extended pairing}
\UU^{\geq} \otimes \UU^{\leq} \xrightarrow{\langle \cdot , \cdot \rangle} \BF
\end{equation}
which satisfies the following properties: \medskip

\begin{enumerate}[leftmargin=*]

\item[i)] for all $i,j \in I$, we have:
$$
\Big \langle h^+_i(z), h^-_j(w) \Big \rangle = \frac {\zeta_{ij} \left(\frac zw \right)}{\zeta_{ji} \left(\frac wz \right)} 
$$
where the right-hand side is expanded as $|z| \gg |w|$, and:

\medskip

\item[ii)] for all $i,j \in I$ and $d,k \in \BZ$, we have: 
$$
\Big \langle e_{i,d}, f_{j,k} \Big \rangle = \delta_j^i \delta_{d+k}^0
$$
This property implies that the restriction of  $\langle \cdot, \cdot \rangle$ to $\UUp \otimes \UUm$ coincides with \eqref{eqn:pairing shuffle} under the identifications $\UUp \cong \CS$ and $ \UUm\cong \CS^{\op}$.

\medskip

\end{enumerate}

\noindent Summarizing the discussion above, we may define by the usual Drinfeld double procedure (see \cite[Subsection 4.15]{N} for a review) the following algebra. \medskip

\begin{definition}
\label{def:extendeddouble}
The \textbf{generic quantum loop group} $\mathbf{U}$ is the $\mathbb{F}$-algebra generated by elements:
$$
\{e_{i,k}, f_{i,k}, h^\pm_{i,l}\; | i \in I, k \in \BZ, l \geq 0\}
$$
satisfying the fact that $h_{i,0}^\pm$ are invertible, as well as the quadratic relations: 
\begin{align*}
&e_i(z) e_j (w) \zeta_{ji} \left( \frac wz \right) = e_j(w) e_i(z) \zeta_{ij} \left(\frac zw \right) \\ 
&f_j(w) f_i(z) \zeta_{ji} \left(\frac wz \right)  = f_i(z) f_j (w) \zeta_{ij} \left( \frac zw \right) 
\end{align*}
for all $i,j \in I$, the cubic relations:
\begin{multline}
\frac {\tzeta_{ii}\left(\frac {x_2}{x_1} \right) \tzeta_{ji} \left( \frac {y}{x_1} \right) \tzeta_{ji} \left( \frac {y}{x_2} \right)}{\left(1- \frac {x_2}{x_1q}\right) \left(1-\frac {yq}{x_2t_{\edge}} \right)} \cdot e_i(x_1) e_i(x_2) e_j(y)  \\
+ \frac {\tzeta_{ii}\left(\frac {x_1}{x_2} \right) \tzeta_{ji} \left( \frac y{x_2} \right) \tzeta_{ij} \left( \frac {x_1}y \right) \left(- \frac {x_2 t_{\edge}}{y} \right) \left(- \frac y{x_1} \right)^{\delta_j^i}}{\left(1- \frac {yq}{x_2t_{\edge}}\right) \left(1-\frac {x_1 t_{\edge}}{y} \right)} \cdot e_i(x_2) e_j(y) e_i(x_1)  \\
+ \frac {\tzeta_{ii}\left(\frac {x_2}{x_1} \right) \tzeta_{ij} \left( \frac {x_1}y \right) \tzeta_{ij} \left( \frac {x_2}y \right) \left( \frac {x_2t_{\edge}}{yq} \right)\left( \frac {y^2}{x_1x_2} \right)^{\delta_j^i}}{\left(1- \frac {x_2}{x_1 q}\right) \left(1-\frac {x_1 t_{\edge}}{y} \right)} \cdot e_j(y) e_i(x_1) e_i(x_2) = 0
\end{multline}
\begin{multline}
\frac {\tzeta_{ii}\left(\frac {x_2}{x_1} \right) \tzeta_{ji} \left( \frac {y}{x_1} \right) \tzeta_{ji} \left( \frac {y}{x_2} \right)}{\left(1- \frac {x_2}{x_1q}\right) \left(1-\frac {yq}{x_2t_{\edge}} \right)} \cdot  f_j(y)f_i(x_2)f_i(x_1)\\
+ \frac {\tzeta_{ii}\left(\frac {x_1}{x_2} \right) \tzeta_{ji} \left( \frac y{x_2} \right) \tzeta_{ij} \left( \frac {x_1}y \right) \left(- \frac {x_2 t_{\edge}}{y} \right) \left(- \frac y{x_1} \right)^{\delta_j^i}}{\left(1- \frac {yq}{x_2t_{\edge}}\right) \left(1-\frac {x_1 t_{\edge}}{y} \right)} \cdot  f_i(x_1)f_j(y)f_i(x_2) \\
+ \frac {\tzeta_{ii}\left(\frac {x_2}{x_1} \right) \tzeta_{ij} \left( \frac {x_1}y \right) \tzeta_{ij} \left( \frac {x_2}y \right) \left( \frac {x_2t_{\edge}}{yq} \right)\left( \frac {y^2}{x_1x_2} \right)^{\delta_j^i}}{\left(1- \frac {x_2}{x_1 q}\right) \left(1-\frac {x_1 t_{\edge}}{y} \right)} \cdot f_i(x_2)f_i(x_1)f_j(y) = 0
\end{multline}
for any edge $\overline{E} \ni e = \oij$, as well as:
\begin{align}
&e_i(z) h^\pm_j(w) = h^\pm_j(w) e_i(z) \frac {\zeta_{ij} \left( \frac z{w} \right)}{\zeta_{ji} \left( \frac {w}z \right)}  \qquad (\text{for\;}|z^{\pm 1}| \ll |w^{\pm 1}|)\\
&f_i(z) h^\pm_j(w) = h^\pm_j(w) f_i(z) \frac {\zeta_{ji} \left( \frac {w}z \right)}{\zeta_{ij} \left( \frac z{w} \right)} \qquad (\text{for\;}|z^{\pm 1}| \ll |w^{\pm 1}|)
\end{align} 
and:
$$
[e_{i,d}, f_{j,k}] = \delta_{j}^i \cdot \begin{cases} - h^+_{i,d+k} &\text{if } d+k > 0 \\[4pt] h_{i,0}^{-} - h_{i,0}^{+} &\text{if } d+k = 0 \\[4pt] h^-_{i,-d-k} &\text{if } d+k<0 \end{cases}
$$
for all $i,j \in I$. Formulas \eqref{eqn:cop 1}, \eqref{eqn:cop 2}, \eqref{eqn:cop 3} endow $\mathbf{U}$ with the structure of a bialgebra. Finally, there is a triangular decomposition as $\BF$-vector spaces: 
$$
\mathbf{U} \simeq \mathbf{U}^+ \otimes \mathbf{U}^0\otimes \mathbf{U}^-
$$
where $\mathbf{U}^+, \mathbf{U}^0, \mathbf{U}^-$ are respectively generated by $\{e_{i,k}\}_{i,k}$, $\{h_{i,l}^\pm\}_{i,l}$ and $\{f_{i,k}\}_{i,k}$.
\end{definition}

\medskip

\noindent As is standard in the theory of quantum groups, one can write down antipode maps that make $\UU^{\geq}$ and $\UU^{\leq}$ into Hopf algebras, and also define central extensions by making the series $h^+(z)$ and $h^-(w)$ ``almost" commute with each other (see \cite{R-matrix} for a survey in the case of the Jordan quiver). We will not describe these in detail. \medskip

\begin{remark}

Since it arises as a Drinfeld double, $\UU$ has an interesting universal $R$-matrix (up to the usual issues involving completions and central extensions that are required to rigorously define it). This object is studied in \cite{N2}, where it is conjectured that it matches the $R$-matrices defined via Nakajima quiver varieties by Aganagic, Maulik, Okounkov and Smirnov (cf.\ \cite{MO_Quantum, OS_Quantum, O_Quantum, AO_Envelopes}). 

\end{remark}

\bigskip

\section{Special values of the parameters}
\label{sec:specialize}

\medskip

\subsection{} An important question is to determine what happens to $\UUp$ and $\CS$ when the parameters $q$ and $\{t_e\}_{e \in E}$ are no longer generic, i.e. when we work over a field different from \eqref{eqn:ground field}. In the present Section, we will provide an answer under Assumption \begin{otherlanguage*}{russian}Ъ\end{otherlanguage*} of \cite{N}, given as follows: 

\medskip

\begin{definition}
\label{def:assumption}

Let $\BK$ be a field endowed with elements $q$ and $\{t_e\}_{e \in E}$, such that there exists a field homomorphism:
$$
\rho : \BK \longrightarrow \BC
$$
for which $|\rho(q)| < |\rho(t_e)| < 1$ for all $e \in E$. \\

\end{definition}

\noindent We emphasize the fact that $\BK$ will henceforth refer to a choice of a field \textit{together} with the elements $q,t_e$ as above. With this in mind, we may define:
$$
{_\BK\CV} \quad \text{and} \quad {_\BK\tUUp}
$$
simply by replacing $\BF$ with $\BK$ in the definition of $\CV$ and $\tUUp$ of  Definitions \ref{def:shuf big} and \ref{def:big quantum}, respectively. All the results of Section \ref{sec:big} will continue to hold in the present context, notably the existence of a pairing:
\begin{equation}
\label{eqn:restricted pairing}
{_\BK\tUUp} \otimes {_\BK\CV^{\op}} \xrightarrow{\langle \cdot , \cdot \rangle} \BK
\end{equation}
analogous to \eqref{eqn:pairing}. This pairing is nondegenerate in the second variable by arguments similar to those in Section~\ref{sec:big}. It will follow from Theorem~\ref{thm:main restricted} that this pairing is nondegenerate in the first variable as well (we will not use this fact in the proof of Theorem~\ref{thm:main restricted}).

\medskip

\subsection{} However, as soon as we reach Definition \ref{def:shuffle}, we notice that \eqref{eqn:wheel} cannot be the correct definition of the small shuffle algebra anymore (for one thing, it might be the case that $t_e = t_{e'}$ for two distinct edges $e,e'$ between vertices $i$ and $j$, in which case two wheel conditions \eqref{eqn:wheel} actually impose the same restriction on $R$). As shown in \cite{N}, the way to fix this in the context of Definition \ref{def:assumption}, is to set:
$$
{_\BK\CV} \supset {_\BK\CS} = \Big\{R( \dots, z_{i1}, \dots, z_{in_i}, \dots) \text{ such that } \forall i \in I
$$
\begin{equation}
\label{eqn:restricted wheel}
R\Big|_{z_{i2} = qz_{i1}} \quad \text {is divisible by} \quad \prod_{(j,b) \notin \{(i,1), (i,2)\}} \prod_{\overline{E} \ni e = \oij} (z_{jb} - t_{e} z_{i1}) \Big\}
\end{equation}
It is easy to see that \eqref{eqn:restricted wheel} is equivalent to the following condition for all $i,j \in I$ and all $\gamma \in \BK$ (below, we replace $z_{j1}$ by $z_{i3}$ if $i = j$):
\begin{equation}
\label{eqn:restricted wheel equivalent}
R\Big|_{z_{i2} = qz_{i1}} \quad \text {is divisible by} \quad (z_{j1} - \gamma z_{i1})^{\flat_{ij}(\gamma)}
\end{equation}
where $\flat_{ij}(\gamma)$ denotes the number of edges $e = \oij$ in $\overline{E}$ for which $t_e = \gamma$. This is why if the parameters $\{t_e\}_{e \in \overline{E}}$ are all distinct, we recover the wheel conditions \eqref{eqn:wheel}.

\medskip

\noindent As shown in \cite{N}, the definition \eqref{eqn:restricted wheel} ensures that ${_\BK\CS}$ is generated by $\{z_{i1}^d\}_{i \in I}^{d\in \BZ}$, and that there is a non-degenerate pairing:
\begin{equation}
\label{eqn:restricted pairing shuffle}
{_\BK\CS} \otimes {_\BK\CS}^{\op} \longrightarrow \BK
\end{equation} 
analogous to \eqref{eqn:pairing shuffle}.

\medskip

\subsection{} We will now define a quotient:
$$
{_\BK\tUUp} \twoheadrightarrow {_\BK\UUp}
$$
for which the analogue of Theorem \ref{thm:main} holds. Let us recall the notation of Subsection \ref{sub:cubic}, and define for all $i,j \in I$, $0 \leq k \leq \flat_{ij}(\gamma)$, all $\gamma \in \BK$ and all integers $(a,b,c) \in \BZ^3$:
\begin{equation}
\label{eqn:defA restricted}
{_\BK}A^{(i,j,\gamma|k)}_{a,b,c} =
\end{equation}
$$
\left[ \frac{x_1^a \left( \frac{x_2}{q}\right)^b \left(\frac{y}{\gamma}\right)^c }{(1-q)\left(1 - \gamma\right)^{\delta_j^i}\left(1 - \frac {\gamma}q\right)^{\delta_j^i}} \cdot {_{\BK}X}^{(i,j,\gamma|k)}(x_1,x_2,y) \right]_{\text{ct}} \in \ {_\BK\tUU}^+_{2\bs_i+\bs_j,a+b+c} 
$$
where:
\begin{equation}
\label{eqn:defX restricted}
{_\BK}X^{(i,j,\gamma|k)}(x_1,x_2,y) 
\end{equation}
is defined as the formal series:
\begin{align*}
& (-1)^{k-1}\frac {\tzeta_{ii}\left(\frac {x_2}{x_1} \right) \tzeta_{ji} \left( \frac {y}{x_1} \right) \tzeta_{ji} \left( \frac {y}{x_2} \right) \left(\frac {y}{x_1\gamma} \right)}{\left(1- \frac {x_2}{x_1q}\right) \left(1-\frac {yq}{x_2\gamma} \right)^k} \cdot e_i(x_1) e_i(x_2) e_j(y) \\
&+ \sum^{k', k'' \in \BZ_{>0}}_{k'+k'' = k+1}(-1)^{k'-1} \frac {\tzeta_{ii}\left(\frac {x_1}{x_2} \right) \tzeta_{ji} \left( \frac y{x_2} \right) \tzeta_{ij} \left( \frac {x_1}y \right) \left(\frac {x_1 \gamma}{y} \right)^{k''-1} \left( - \frac {yq}{x_1\gamma} \right)\left(- \frac y{x_1} \right)^{\delta_j^i} }{\left(1- \frac {yq}{x_2\gamma}\right)^{k'} \left(1-\frac {x_1 \gamma}{y} \right)^{k''}} \cdot e_i(x_2) e_j(y) e_i(x_1) \\
&+ \frac {\tzeta_{ii}\left(\frac {x_2}{x_1} \right) \tzeta_{ij} \left( \frac {x_1}y \right) \tzeta_{ij} \left( \frac {x_2}y \right) \left( \frac {x_1\gamma}{y} \right)^{k-1} \left(\frac {y^2}{x_1x_2} \right)^{\delta_j^i}}{\left(1- \frac {x_2}{x_1 q}\right) \left(1-\frac {x_1 \gamma}{y} \right)^k} \cdot e_j(y) e_i(x_1) e_i(x_2)  
\end{align*}
If $k \leq \flat_{ij}(\gamma)$, the linear factors in $x_1,x_2,y$ in the denominators above are all canceled by similar linear factors in the numerator, hence \eqref{eqn:defX restricted} is a Laurent polynomial in $x_1,x_2,y$ times the product of the series $e_i(x_1)$, $e_i(x_2)$, $e_j(y)$. Note that when $k = \flat_{ij}(\gamma) = 1$, the expression \eqref{eqn:defA restricted} is slightly different from \eqref{eqn:defA}, but they both yield equivalent sets of cubic relations modulo the ideal of quadratic relations \eqref{eqn:rel quad}. As such, it is straightforward to generalize Proposition \ref{prop:wheel pairing} to the following result. 

\medskip

\begin{proposition}
\label{prop:wheel pairing restricted}

For all $i,j \in I$, $R \in \CV^{\emph{op}}_{2\bs_i+\bs_j}$, $0 \leq k \leq \flat_{ij}(\gamma)$, $\gamma \in \BK$ and $(a,b,c) \in \BZ^3$, we have:
\begin{equation}
\label{eqn:wheel pair restricted}
\Big \langle {_\BK}A^{(i,j,\gamma|k)}_{a,b,c}, R \Big \rangle = \int   \frac {\partial^{k-1} F}{\partial z_{j1}^{k-1}}  \Big|_{z_{j1} = z_{i1} \gamma} \frac {(z_{i1}\gamma)^{k-1}}{(k-1)!} Dz_{i1}
\end{equation}
%(the integral refers to the difference of the residues in the variable $z_{i1}$ at 0 and $\infty$ of the integrand)
where:
$$
F(z_{i1}, z_{j1}) = z_{i1}^{a+b} \left(\frac {z_{j1}}{\gamma}\right)^c R(z_{i1}, z_{i1}q, z_{j1}) \left(1 - \frac {z_{j1}}{z_{i1}}\right)^{\delta_j^i} \left(1 - \frac {z_{j1}}{z_{i1} q}\right)^{\delta_j^i}
$$
(we replace $z_{j1}$ by $z_{i3}$ in the formulas above if $i = j$). 

\medskip

\end{proposition} 

\begin{proof} Consider the following formal series, for all $k \in \BN$:
$$
\delta^{(k)}(x) = \sum_{d \in \BZ} x^{d} {-d \choose k}
$$
(note that $\delta^{(0)} = \delta$). The following property is well-known:
\begin{equation}
\label{eqn:identity}
\left[ \delta \left(\frac {z_{i2}}{z_{i1}q} \right) \delta^{(k)} \left(\frac {z_{j1}}{z_{i1}\gamma} \right) P(z_{i1},z_{i2},z_{j1}) \right]_{\text{ct}} = \int \frac {\partial^{k} P(z_{i1},qz_{i1},z_{j1})}{\partial z_{j1}^{k}} \Big|_{z_{j1} = z_{i1} \gamma} \frac {(z_{i1}\gamma)^{k}}{k!} Dz_{i1}
\end{equation}
for all Laurent polynomials $P(z_{i1}, z_{i2}, z_{j1})$, where $[\dots]_{\text{ct}}$ refers to the constant term in the variables $z_{i1}, z_{i2}, z_{j1}$. With this in mind, formula \eqref{eqn:wheel pair restricted} reduces to the following identity of formal power series (whose $k=1$ case is a closely related version of \eqref{eqn:delta functions}):
\begin{align}
\delta \left(\frac {z_{i2}}{z_{i1}q} \right)\delta^{(k-1)} \left(\frac {z_{j1}}{z_{i1}\gamma} \right) &= ev_{ |z_{i1}| \gg |z_{i2}| \gg |z_{j1}| } \left[\frac {(-1)^{k-1} \left(\frac {z_{i2}}{z_{i1}q} \right) \left(\frac {z_{j1}q}{z_{i2}\gamma} \right)}{\left(1- \frac {z_{i2}}{z_{i1}q}\right) \left(1-\frac {z_{j1}q}{z_{i2}\gamma} \right)^k} \right] \\
&+  ev_{ |z_{i2}| \gg |z_{j1}| \gg |z_{i1}| } \left[ \sum^{k', k'' \in \BZ_{>0}}_{k'+k'' = k+1} \frac {(-1)^{k'-1} \left(\frac {z_{j1}q}{z_{i2}\gamma} \right) \left(\frac {z_{i1}\gamma}{z_{j1}} \right)^{k''-1}}{\left(1- \frac {z_{j1}q}{z_{i2}\gamma}\right)^{k'} \left(1-\frac {z_{i1}\gamma}{z_{j1}} \right)^{k''}}\right]  \\
&+  ev_{ |z_{j1}| \gg |z_{i1}| \gg |z_{i2}| } \left[\frac {\left(\frac {z_{i1}\gamma}{z_{j1}} \right)^{k-1}}{\left(1- \frac {z_{i2}}{z_{i1}q}\right) \left(1-\frac {z_{i1}\gamma}{z_{j1}} \right)^k}\right] \label{eqn:delta functions restricted}
\end{align}
If we relabel $z_{i1}=x$, $z_{j1} = y\gamma$, $z_{i2} = zq$, then formula \eqref{eqn:delta functions restricted} has left-hand side:
\begin{equation}
\label{eqn:lhs}
\sum_{d_1,d_2 \in \BZ} x^{-d_1-d_2} y^{d_1} z^{d_2} {-d_1 \choose k-1}
\end{equation}
Because of the elementary identity
\begin{equation}\label{eqn:taylorexpansion}
ev_{|v| \gg |u|} \left[ \frac 1{\left(1-\frac uv \right)^k} \right] = \sum_{d = 0}^{\infty} \frac {(-u)^d}{v^{d}} {-k \choose d} = (-1)^{k-1} \sum_{d = 0}^{\infty} \frac {u^d}{v^{d}} {-d-1 \choose k-1}
\end{equation}
the right-hand side of \eqref{eqn:delta functions restricted} is equal to:
\begin{align}
&\sum_{d_1 = 1}^{\infty} \sum_{d_2 = 1-d_1}^{\infty} x^{-d_1-d_2} y^{d_1} z^{d_2} {-d_1 \choose k-1} \\ + &\sum^{k', k'' \in \BZ_{>0}}_{k'+k'' = k+1} \sum_{d_2=-\infty}^{-1} \sum_{d_1 = -\infty}^{-d_2-1} x^{-d_1-d_2} y^{d_1} z^{d_2} {d_2 \choose k'-1}{-d_1-d_2 \choose k''-1} \\ + &\sum_{d_1 = -\infty}^{0} \sum_{d_2 = 0}^{\infty} x^{-d_1-d_2} y^{d_1} z^{d_2} {-d_1 \choose k-1} \label{eqn:rhs}
\end{align}
To obtain \eqref{eqn:rhs}, we relabeled certain indices and used the Taylor series expansion \eqref{eqn:taylorexpansion} for $u$ in place of $u/v$.
%$$
%\frac 1{(1-u)^k} = \sum_{d=0}^{\infty} (-u)^d {-k \choose d} = (-1)^{k-1} \sum_{d=0}^{\infty} u^d {-d-1 \choose k-1}
%$$
Using a simple combinatorial identity involving binomial coefficients, expressions \eqref{eqn:lhs} and \eqref{eqn:rhs} are easily seen to be equal to each other, thus establishing \eqref{eqn:delta functions restricted}. 
\end{proof}

\subsection{} Formula \eqref{eqn:wheel pair restricted} shows that \eqref{eqn:restricted wheel equivalent} holds for $R$ if and only if:
\begin{equation}
\label{eqn:pairing vanish restricted}
\Big \langle {_\BK}A^{(i,j,\gamma|k)}_{a,b,c}, R \Big \rangle = 0 \qquad \forall k \leq \flat_{ij}(\gamma) \text{ and } a,b,c \in \BZ
\end{equation}
Therefore, generalizing \eqref{eqn:small quantum}, we define:
\begin{equation}
\label{eqn:quotient restricted}
{_\BK\UUp} = {_\BK\tUUp} \Big / {_\BK}J
\end{equation}
where ${_\BK}J$ denotes the two-sided ideal of the quadratic quantum loop group ${_\BK\tUUp}$ generated by the elements \eqref{eqn:defA restricted}, for all $i,j \in I$, $\gamma \in \BK$,  $k \in \{1,\dots,\flat_{ij}(\gamma)\}$ and $a,b,c \in \BZ$. 

\medskip

\begin{remark}
\label{rem:consider}

One could (and should) define ${_\BK}J$ in \eqref{eqn:quotient restricted} to be the two-sided ideal generated by finitely many of the elements \eqref{eqn:defA restricted} in every fixed degree. Indeed, we simply need to use enough $a,b,c$ in \eqref{eqn:pairing vanish restricted} to ensure that $R$ satisfies property \eqref{eqn:restricted wheel equivalent}. The latter property is equivalent \footnote{We are not saying anything deep here, simply the fact that for all constants $\alpha,\beta,\gamma$, we have $f(z\alpha,z\beta,z\gamma) = 0$ if and only if
$$
\int \left( z_1^a z_2^b z_3^c f(z_1,z_2,z_3) \Big|_{z_1 = z\alpha, z_2 = z\beta, z_3 = z\gamma} \right) Dz = 0
$$
for a single triple $(a,b,c)$ of any fixed $d  = a+b+c \in \BZ$.} to the vanishing of the r.h.s. of equation \eqref{eqn:wheel pair restricted} (for all applicable $i,j,k,\gamma$) for a single triple $(a,b,c)$ of any fixed $d = a+b+c \in \BZ$. Thus, in \eqref{eqn:pairing vanish restricted}, it is enough to consider (for all applicable $i,j,k,\gamma$) a single triple $(a,b,c)$ of any fixed $d = a+b+c \in \BZ$. Compare with Proposition \ref{prop:wheel pairing}.

\end{remark}

\medskip

\noindent The analogue of Theorem~\ref{thm:main} in the setting at hand is the following.

\medskip

\begin{theorem}
\label{thm:main restricted}

The assignment $\te_{i,d} \mapsto z_{i1}^d$ for all $i \in I, d \in \BZ$ induces an algebra isomorphism: 
$$
{_\BK\UUp} \xrightarrow{\sim} {_\BK\CS}.
$$

\end{theorem}

\medskip

\noindent We will now sketch the main points of the proof of Theorem~\ref{thm:main restricted}, leaving the details to the interested reader. Just like the proof of Theorem~\ref{thm:main} started with the non-degenerate pairing \eqref{eqn:pairing shuffle}, in the case at hand we start from the pairing \eqref{eqn:restricted pairing shuffle}, whose non-degeneracy was the main reason for introducing the specific conditions \eqref{eqn:restricted wheel equivalent} in \cite{N}. Then the proof of Theorem~\ref{thm:main} carries through as stated, once one replaces the elements \eqref{eqn:defA} of the quadratic quantum loop group by the elements \eqref{eqn:defA restricted} (in both situations, one needs only consider finitely many such elements in any given degree; see Remark~\ref{rem:consider}).

\bigskip

\section{The spherical Hall algebra of a generic curve}\label{sec:genusg}

\medskip

\subsection{}

Let us fix $g \in \BN$ and specialize the quiver to $Q=S_g$, the quiver with one vertex and $g$ loops. In this special case, we will show how to connect the extended quantum loop group $\UU^{\geq}$ of Definition \ref{def:extended} with the Hall algebra of a smooth projective curve of genus $g$ over a finite field. In the present Section, we will treat the case of a generic curve (in the sense of having distinct Weil numbers), and then in the next Section we will show how to adapt to the case of an arbitrary curve.

\medskip

\noindent In more detail, let $X$ be a smooth, geometrically connected projective curve of genus $g$ defined over the finite field $\mathbf{k}=\BF_{q^{-1}}$ (here the unusual choice of notation for the cardinality of the finite field is made for the purpose of compatibility with \eqref{eqn:def zeta}). Let $Coh_{r,d}(\finite)$ denote the groupoid of coherent sheaves on $X$ of rank $r$ and degree $d$. We refer to 
\cite[Lectures~1, 4]{SLectures} for the definition of the Hall algebra of the category of coherent sheaves on $X$. As a vector space:
\begin{equation}
\label{eqn:def hall}
\mathbf{H}_X=\left(\bigoplus_{(r,d) \in (\BZ^2)^+} \text{Fun}_0(Coh_{r,d}(\finite), \BQ)\right) \otimes \BC[\kappa^{\pm 1}] 
\end{equation}
is the space of finitely supported functions on $Coh(\finite)=\bigsqcup_{r,d} Coh_{r,d}(\finite)$, where:
$$
(\BZ^2)^+:=\{(r,d) \in \BZ^2\;|\;r \geq 0, d \in \BZ \text{ such that } d \geq 0 \text{\;if\;} r=0\}.
$$
The element $\kappa$ is usually denoted by $\kappa_{1,0}$ in the literature, and satisfies the following commutation relations :
$$\kappa f \kappa^{-1}=q^{r(g-1)}f, \qquad f \in \text{Fun}_0(Coh_{r,d}(\finite), \BQ);$$
 one often adds a central element (denoted by $\kappa_{0,1}$) to $\mathbf{H}_X$, although we will not do so, as it only plays a significant part when considering the double Hall algebra. 

\medskip

\begin{definition}
\label{def:spherical}

The spherical Hall algebra $\mathbf{H}^{sph}_X$ is the subalgebra of $\mathbf{H}_X$ generated by $\kappa^{\pm 1}$ together with the following elements:
$$
1_{1,d}^{vec}=\chi_{Pic_d(\finite)}\ , \qquad 1_{0,l}=\chi_{Coh_{0,l}(\finite)}\ , \qquad (d \in \BZ, l \geq 1)
$$
where $\chi_Y$ stands for the characteristic function of a subgroupoid $Y \subset Coh(\finite)$. 

\end{definition}

\medskip

\noindent Let us consider the subalgebras:
\begin{align*} 
&\mathbf{H}_X^{sph} \supset \mathbf{H}^{sph,+}_X = \BC \langle  1^{vec}_{1,d}\;|\; d \in \BZ \rangle \\
&\mathbf{H}_X^{sph} \supset \mathbf{H}^{sph,0}_X = \BC\langle \kappa^{\pm 1}, 1_{0,l}\;|\; l \geq 1 \rangle 
\end{align*}
In Proposition \ref{prop:semidirect}, we will recall the way to realize $\mathbf{H}_X^{sph}$ as a semidirect product of $\mathbf{H}^{sph,+}_X$ with $\mathbf{H}^{sph, 0}_X$ . In order to set this up, it is convenient to introduce new generators $\{T_{0,l}\}, \{\theta_{0,l}\}, l \geq 1$ of $\mathbf{H}^{sph,0}_X$  through the relations:
\begin{align*}
&1+\sum_{l = 1}^{\infty} 1_{0,l}s^l = \exp\left( \sum_{l = 1}^{\infty} T_{0,l} \frac{s^l}{[l]}\right) \\
&1+\sum_{l = 1}^{\infty} \theta_{0,l}s^l = \exp\left(\sum_{l = 1}^{\infty}(q^{-1/2}-q^{1/2})T_{0,l}s^l\right)
\end{align*}
where $[l]=(q^{l/2}-q^{-l/2})/(q^{1/2}-q^{-1/2})$. It is known that $\mathbf{H}_X^{sph,0}$ is a free commutative polynomial algebra in $\kappa^{\pm 1}$ and any one of the families of generators $\{1_{0,l}\}_{l \geq 1}$, $\{T_{0,l}\}_{l \geq 1}$ or $\{\theta_{0,l}\}_{l \geq 1}$.

\medskip

\subsection{}

Let $\sigma_1, \overline{\sigma_1}, \ldots, \sigma_g, \overline{\sigma_g}$ denote the Weil numbers of $X$, paired up such that $\sigma_e \overline{\sigma_e} = q^{-1}$ for all $e = 1, \dots, g$. Therefore, we may write:
\begin{equation}
\label{eqn:weil numbers}
t_e = \frac 1{\sigma_e} \qquad \text{and} \qquad t_{e^*} = \frac 1{\overline{\sigma_{e}}}
\end{equation}
for all $e = 1, \dots, g$, and this would be compatible with \eqref{eqn:inverse t}. With this notation, the particular case of the rational function \eqref{eqn:def zeta} when $Q$ is the $g$ loop quiver (one vertex with $g$ loops) is a renormalized form of the zeta function of $X$:
\begin{equation}
\label{eqn:zeta of curve}
\zeta_X(x)=\frac{1-xq^{-1}}{1-x}\prod_{e=1}^g (\sigma_e-x)(1-\overline{\sigma_e}x^{-1})
\end{equation}
(explicitly, \eqref{eqn:zeta of curve} is equal to \cite[(1.22)]{SV} times $(-1)^{g-1}$). Finally, set:
$$
E(z)=\sum_{d \in \BZ} 1^{vec}_{1,d}z^{-d}, \qquad H^+(z)=\kappa \left( 1 + \sum_{l =  1}^{\infty} \theta_{0,l}z^{-l}\right)
$$

\medskip

\begin{proposition} 
\label{prop:semidirect}

The multiplication map $\mathbf{H}^{sph,+}_X \otimes \mathbf{H}^{sph,0}_X \to \mathbf{H}^{sph}_X$ is a vector space isomorphism.  In addition, as an algebra, $\mathbf{H}^{sph}_X$ is generated by $\mathbf{H}^{sph,+}_X$ and $\mathbf{H}^{sph,0}_X$ modulo the following relations:
\begin{equation}\label{eqn:EH}
E(z)H^+(w)=H^+(w)E(z) \frac{\zeta_X\left(\frac{z}{w}\right)}{\zeta_X\left(\frac{w}{z}\right)} \qquad (\text{for\;}|w| \gg |z|)
\end{equation}
\end{proposition}
\begin{proof} See \cite[Corollary~1.4 and Section~1.5]{SV}. Note that $1_{1,d}^{vec} \kappa = q^{1-g} \kappa 1_{1,d}^{vec}$.
\end{proof}

\medskip

\noindent
Concerning the structure of $\mathbf{H}^{sph,+}_X$, we have the following result. Let $\CV_X$ be the $\BC$-algebra defined by the relations \eqref{eqn:shuf prod} when $Q$ is the $g$-loop quiver (and the equivariant parameters $t_e$ are specialized as in \eqref{eqn:weil numbers}), and let $\CV_X^{sph}$ be the subalgebra of $\CV_X$ generated by its horizontal degree 1 pieces $\{\CV_{X|1,d}\}_{d \in \BZ}$.

\medskip

\begin{theorem}[{\cite[Theorem~1]{SV}}]\label{thm:g1}The assignment $1^{vec}_{1,d} \mapsto z_1^d \in \CV_{X|1,d}$ for $d \in \BZ$ extends to an algebra isomorphism $\mathbf{H}^{sph,+}_X \xrightarrow{\sim} \CV^{sph}_X$.
\end{theorem}

\medskip

\subsection{} As long as the Weil numbers of $X$ are distinct (which implies that $\{t_e\}_{e \in \overline{E}}$ are distinct complex numbers), the small shuffle algebra behaves just like in the case of generic parameters. In the language of Section \ref{sec:specialize}, this is because conditions \eqref{eqn:restricted wheel} are still none other than the usual 3-variable wheel conditions \eqref{eqn:wheel}. Thus, the following is proved just like \cite[Theorem~1.1]{N}.

%We now pass to the formal setup. Let $\CV_g$ be the big shuffle algebra (over $\mathbb{F}$) associated to the quiver with one vertex and $g$ loops. Let $\CV^{sph}_g$ be the subalgebra of $\CV_g$ generated by $\bigoplus_d \CV_{g|1,d}$. Let us put $\mathbf{H}^{sph,0}_g=\mathbb{F}[\kappa^{\pm 1}, \theta_{0,l}\;; l\geq 1]$. We define the \textbf{generic spherical Hall algebra of genus $g$} as the semi-direct product 
%$$
%\mathbf{H}^{sph}_g=\CV_g^{sph} \rtimes \mathbf{H}_g^{sph,0}
%$$
%using relation \eqref{eqn:EH}, where, slightly abusing notation, we have now set 
%$$
%1^{vec}_d\coloneqq z^d \in \CV_{g|1,d}
%$$ 
%Let $\CS_g$ be the small shuffle algebra, i.e the subalgebra of $\CV_g$ determined by the 3-wheel conditions \eqref{eqn:wheel}. The crucial result here is the following

\medskip

\begin{theorem} \label{thm:g2} We have $\CV^{sph}_X = \CS_X$, namely the subalgebra of $\CV_X$ determined by the 3-variable wheel conditions \eqref{eqn:wheel} (but with $t_e$ specialized as in \eqref{eqn:weil numbers}). \medskip

\end{theorem}

\begin{proof}[Proof of Theorem~\ref{thm:genusg}] By combining Theorem~\ref{thm:g1}, Theorem~\ref{thm:g2} and Theorem~\ref{thm:intro}, we only need to check that the collections of relations \eqref{eqn:serre a} and \eqref{eqn:genusg4} are equivalent (assuming relations \eqref{eqn:genusg1}, \eqref{eqn:genusg2}, \eqref{eqn:genusg3} hold). By \eqref{eqn:genusg2}, taking commutators with the symbols $\theta_{0,l}$ means that formula \eqref{eqn:genusg4} is equivalent to:
$$
\Big[p(x,y,z) (x+z)(xz-y^2)Q_e(x,y,z) E(x)E(y)E(z) \Big]_{\text{ct}} = 0
$$
for any symmetric Laurent polynomial $p$. Thus, we need to prove that for any given $R(z_1,z_2,z_3) \in \CV_X$, the wheel conditions \eqref{eqn:wheel} are equivalent to the relations:
$$
\left \langle \Big[P(x,y,z) (x+z) \left(\frac 1y - \frac y{xz} \right) Q_e(x,y,z) E(x)E(y)E(z) \right]_{\text{ct}} , R \Big \rangle = 0
$$
for any symmetric Laurent polynomial $P(x,y,z) = p(x,y,z)xyz$. By definition of the pairing \eqref{eqn:pairing formula}, the condition above reads precisely:
\begin{multline}
\label{eqn:twolines}
\int_{|z_1| \gg |z_2| \gg |z_3|} \frac {Q_e \left(z_1,z_2,z_3 \right)}{\zeta_X \left( \frac {z_2}{z_1} \right) \zeta_X \left( \frac {z_3}{z_1} \right) \zeta_X \left( \frac {z_3}{z_2} \right)} \\
(z_1 + z_3)\left(\frac 1{z_2} - \frac {z_2}{z_1z_3} \right) (PR)(z_1,z_2,z_3) Dz_1Dz_2Dz_3 = 0
\end{multline}
By the definition of $Q_e$, the rational function on the first row above is equal to:
$$
\frac 1{\zeta_1^{(e)} \left( \frac {z_2}{z_1} \right) \zeta_1^{(e)} \left( \frac {z_3}{z_1} \right) \zeta_1^{(e)} \left( \frac {z_3}{z_2} \right)}
$$
where $\zeta_1^{(e)}$ is the rational function \eqref{eqn:def zeta} for the Jordan quiver, associated to the single equivariant parameter $t_e$ (besides $q$). This reduces the problem to the case $g = 1$, in which case the equivalence of \eqref{eqn:twolines} and the wheel conditions \eqref{eqn:wheel} follows from the main results of \cite{S} and \cite{Shuf}. In more detail, it is straightforward to show that for any symmetric Laurent polynomials $P$ and $R$, we have:
$$
\int_{|z_1| \gg |z_2| \gg |z_3|} \frac {(z_1 + z_3)\left(\frac 1{z_2} - \frac {z_2}{z_1z_3} \right) (PR)(z_1,z_2,z_3)}{\zeta_1^{(e)} \left( \frac {z_2}{z_1} \right) \zeta_1^{(e)} \left( \frac {z_3}{z_1} \right) \zeta_1^{(e)} \left( \frac {z_3}{z_2} \right)} Dz_1Dz_2Dz_3 = 
$$
$$
= \int \frac {(PR)(x,xt_e,xq)- (PR)(x,\frac {xq}{t_e},xq)}{\left(\frac q{t_e} - t_e\right)\left(\frac 1q - t_e \right) \left(\frac 1q - \frac q{t_e} \right) q^{-3}} Dx
$$
(it suffices to do so when $PR(z_1,z_2,z_3) = \text{Sym }z_1^az_2^bz_3^c$ for various $a,b,c \in \BZ$, in which case the formula above is a straightforward computation involving power series). Clearly, the vanishing of the right-hand side of the equation above (for all symmetric Laurent polynomials $P$) is equivalent to \eqref{eqn:wheel}. 
\end{proof}

\bigskip

\section{The whole Hall algebra of an arbitrary curve}\label{sec:wholeHall}

\medskip

\subsection{} We now briefly explain how one can adapt and mix the results in the previous two Sections to give a presentation, first of the spherical Hall algebra of an \textit{arbitrary} smooth projective curve of genus $g$, and then of a much larger subalgebra of the Hall algebra of such a curve (under the technical hypothesis that the cotangent bundle $\Omega_X$ admits a square root). Recall that $X$ is a smooth, projective, geometrically connected curve of genus $g$ defined over a finite field $\mathbb{F}_{q^{-1}}$, whose Weil numbers are denoted as in \eqref{eqn:weil numbers} and whose renormalized $\zeta$ function is given in \eqref{eqn:zeta of curve}. 

\medskip

\noindent Let us first consider the spherical Hall algebra of $X$. Observe that because: 
$$
|\sigma_e|=q^{- \frac 12} \quad \Rightarrow \quad |t_e| = q^{\frac 12}
$$
for all $e$, we are in the situation of Assumption \begin{otherlanguage*}{russian}Ъ\end{otherlanguage*} (i.e. Definition \ref{def:assumption}). We can therefore apply the results of Theorem~\ref{thm:main restricted} and obtain a presentation of $\mathbf{H}^{sph}_X$ by generators and relations, the form of which only depends on the multiplicities of the various Weil numbers, i.e. on the potential numerical coincidences between the elements of the multiset $\{t_e, t_{e^*}\}_{e = 1, \dots, g}$. This leads to the following result. We use the same definitions as in Theorem~\ref{thm:genusg} for the generating series $E(z), H(z)$, and we let $X^{(\gamma|k)}(x_1,x_2,x_3)$ be defined as in \eqref{eqn:defX restricted}, but with $\widetilde{\zeta}_{ij}(u)$ replaced by $(1-u)\zeta_X(u)$ and $e_i(u),e_j(u)$ both replaced by $E(u)$.

\medskip

\begin{theorem}\label{thm:genusgarb} Let $\gamma_1, \ldots, \gamma_s$ be the \textit{distinct} elements of $\{t_1, \ldots, t_g, q/t_1, \dots, q/t_g\}$, and let $k_1, \ldots, k_s$ stand for their respective multiplicities. Then $\mathbf{H}_X^{sph}$ is isomorphic to the algebra generated by $\kappa^{\pm 1}, \theta_{0,l},  1^{vec}_{d}$ for $l\geq 1, d \in \BZ$ subject to the relations \eqref{eqn:genusg1}, \eqref{eqn:genusg2}, \eqref{eqn:genusg3} and the relations:
\begin{equation}
\label{eqn:multiplerel}
X^{(\gamma_i|k)}(x_1,x_2,x_3)=0,
\end{equation}
for all $i \in \{1,\dots,s\}$ and $k \in \{1,\dots,k_i\}$. 
\end{theorem}

\medskip

\subsection{}
Let us keep the same notation as above concerning the curve $X$, and consider the \textit{entire} Hall algebra $\mathbf{H}_X$ (i.e. not just the spherical part). For the most part, we will drop the subscript $X$ from the notation of the Hall algebra, as the curve will be fixed. We will recall a few features of $\mathbf{H}$ and refer to \cite{KSV} for precise definitions and more details. We will assume that $X$ has a theta characteristic, i.e that $\Omega_X$ admits a square root $\Omega_X^{1/2}$ (see \cite[Definition 3.9]{KSV}). Let $\mathbf{H}^{+}, \mathbf{H}^{0}$ be the subalgebras of functions \eqref{eqn:def hall} on the stacks of vector bundles and torsion sheaves on $X$, respectively. We will write $\mathbf{H}_{r,d}$, $\mathbf{H}^{+}_{r,d}$, etc, for the graded pieces of the aforementioned Hall algebras, corresponding to sheaves of rank $r$ and degree $d$. There is a canonical decomposition of the commutative algebra $\mathbf{H}^{0}$ according to support: 
$$
\mathbf{H}^{0} = \bigotimes_{x \in X(\finite)} \mathbf{H}_x, \qquad \text{where} \quad \mathbf{H}_x=\mathbb{C}[T_{x,1}, T_{x,2}, \ldots]
$$ 
is generated by primitive elements $T_{x,l} \in \mathbf{H}_{0,l \cdot deg(x)}$, for $l \geq 1$.
Next, there is a left action $\cdot$ of $\mathbf{H}^{0}$ on $\mathbf{H}^{+}$ by \textit{Hecke operators}, which is defined  as the composition:
$$
\mathbf{H}^{0} \otimes \mathbf{H}^{+} \stackrel{m}{\longrightarrow} \mathbf{H} \stackrel{\pi}{\longrightarrow} \mathbf{H}^{+}
$$
where $m$ is the multiplication map and $\pi$ is the natural projection map. It satisfies: 
$$
T_{x,l} \cdot f = [T_{x,l}, f] \qquad \forall\; x \in X(\finite), \ l \geq 1.
$$

\medskip

\subsection{} An element $f \in \mathbf{H}^{+}$ is called cuspidal if the standard coproduct on $\mathbf{H}_X$ (cf. \cite[Theorem 3.3]{K}) satisfies: 
$$
\Delta(f) \in (f \otimes 1 ) \oplus  (\mathbf{H}^{0} \otimes \mathbf{H}^{+}).
$$
We will denote by $\mathbf{H}^{cusp}$ the subspace of cuspidal functions. It is finite-dimensional in any fixed rank $r$ and degree $d$. It is known that $\mathbf{H}^{cusp}$ is a minimal generating subspace of $\mathbf{H}^{+}$, and that $\mathbf{H}^{cusp}$ is stable under the Hecke action of $\mathbf{H}^0$. Cuspidal functions are also closely related to the following standard construction in the theory of automorphic forms. Let $\chi : \mathbf{H}^{0} \to \mathbb{C}$ be an algebra character. A cuspidal eigenform of rank $r$ and of eigenvalue $\chi$ is a nontrivial formal infinite sum $f=\sum_{d} f_d \in \prod_d \mathbf{H}^{+,cusp}_{r,d}$ such that: 
$$
h \cdot f=\chi(h)f, \qquad  \forall\;h \in \mathbf{H}^{0}.
$$ 
There is an action of $\mathbb{C}^*$ on the set of pairs $(\chi,f)$ as above given by:
$$
t \cdot \chi(h)=t^{deg(h)}\chi(h), \qquad t \cdot \sum f_d=\sum t^df_d.
$$
Let $\Sigma_r$ be the set of all eigenvalues of cuspidal eigenforms of rank $r$. The strong multiplicity one theorem says that $\Sigma_r$ is a finite union of $\mathbb{C}^*$-orbits, and that for any $\chi \in \Sigma_r$ there exists a unique (up to scalar) eigenform $f_\chi$ of eigenvalue $\chi$. To a pair of characters $(\chi,\chi') \in \Sigma_r \times \Sigma_s$ one associates a Rankin-Selberg L-function $\LHom(\chi,\chi';z)$ which is known to enjoy the following properties (see \cite[App. B]{Lafforgue} for the formulas and notation below; recall that our base field is $\finite=\mathbb{F}_{q^{-1}}$):

\medskip

\begin{enumerate}[leftmargin=*]
\item[i)] $\LHom(t^{-deg}\chi,\chi';z)=\LHom(\chi,t^{deg}\chi';z)=\LHom(\chi,\chi';tz)$.

\medskip

\item[ii)] If $\mathbb{C}^*\cdot \chi \neq \mathbb{C}^*\cdot \chi'$ then $\LHom(\chi,\chi';z)$ is a polynomial in $z$ of degree $2(g-1)rs$ and constant term $1$, all of whose zeros are of complex norm $q^{-1/2}$.

\medskip

\item[iii)] If $\chi=\chi'$ then there exists a positive integer $d_\chi$ such that:
$$
\LHom(\chi,\chi';z)=\frac{P(z)}{(1-z^{d_\chi})(1-(z/q)^{d_\chi})}
$$
where $P(z)$ is a polynomial in $z$ of degree $2(g-1)rs + 2d_\chi$ and constant term $1$,  all of whose zeros are of complex norm $q^{-1/2}$.
\end{enumerate}
\noindent In i) above, $t^{-\deg}$ refers to the function on $\mathbf{H}^{0}$ given by $h \mapsto t^{-\deg(h)}$. The statements in ii) and iii) concerning the norms of the zeros of $P(z)$ is known as the generalized Riemann hypothesis and is proved by Lafforgue (see \cite[Thm. VI.10]{Lafforgue}).

\medskip

\subsection{} A cuspidal eigenform $f_\chi$ is called \textbf{absolutely cuspidal} if $d_\chi=1$. We will denote by $\mathbf{H}^{abs}$ the subalgebra of $\mathbf{H}$ generated by $\mathbf{H}_{x}$ for $x \in X(\finite)$ and the collection of Fourier modes of the absolutely cuspidal eigenforms. We let $\Sigma_r^{abs} \subset \Sigma_r$ be the collection of eigenvalues of absolutely cuspidal eigenforms.
We will now give a full presentation of $\mathbf{H}^{abs}$ (see Remark~\ref{rem:absHall} for some motivation in considering this specific subalgebra). To this end, for $\chi \in \Sigma_r$ consider the generating series: 
$$
E_\chi(z)=\sum_{d \in \BZ} f_d z^{-d} \in \mathbf{H}^{abs}[[z^{\pm 1}]]
$$
where $f$ is the associated cuspidal eigenform (well-defined up to a scalar). The relevance of Rankin-Selberg L-functions in our context stems from the following: as in \cite[Theorem 3.3]{K}, there exists a series $\Psi_\chi(z) \in \mathbf{H}^0[[z^{-1}]]$ such that:
\begin{equation}\label{eqn:coprodcuspidal}
\Delta(E_\chi(z)) = E_\chi(z) \otimes 1 + \kappa^r \Psi_\chi(z) \otimes E_\chi(z)
\end{equation}  
and for any $\chi, \chi' \in \Sigma:=\bigsqcup_r \Sigma_r$,
\begin{equation}\label{eqn:ratioLfunctions}
E_{\chi}(z) \Psi_{\chi'}(w) = \Psi_{\chi'}(w)E_{\chi}(z) \frac {\LHom(\chi',\chi;qz/w)}{\LHom(\chi',\chi;z/w)}.
\end{equation}
The properties above and the shape of the L-functions $\LHom(\chi,\chi;z)$ for absolutely cuspidal characters $\chi$ suggest an isomorphism between $\mathbf{H}^{abs}$ and the shuffle algebras considered in the present paper. To make this precise, consider the following renormalization of $\LHom$, for a pair of characters $(\chi,\chi') \in \Sigma^{abs}_r \times \Sigma^{abs}_s$~:
\begin{equation}\label{eqn:deflambda}
\zeta_{\chi\chi'}(z) = \begin{cases}
\theta_{\chi,\chi'} z^{(g-1)rs}\LHom(\chi,\chi';z^{-1}) & \text{if\;} \BC^*\cdot \chi \neq \BC^* \cdot \chi'\\[8pt]
\left(1-\frac zq \right) \left(1-\frac 1{zq} \right)z^{(g-1)r^2}\LHom(\chi,\chi;z^{-1}) & \text{if\;} \chi=\chi'
\end{cases}
\end{equation} 
where 
$$
\theta_{\chi,\chi'}= {}^{\pi}(\chi^* \boxtimes \chi')\big(\Omega_X^{1/2}\big)
$$
and $^{\pi}(\chi^* \boxtimes \chi')$ is defined as in \cite[Remark 3.3]{KSV}. Moreover, we have for all $t$:
\begin{equation}
\label{eqn:rescale}
\zeta_{t^{-deg}\chi, \chi'}(z)=\zeta_{\chi,t^{deg}\chi'}(z)=\zeta_{\chi\chi'}(t^{-1}z)
\end{equation}
Using the functional equation for Rankin-Selberg L-functions (\cite[Proposition 3.7]{KSV}), it is straightforward to check that:
$$
\frac{\zeta_{\chi\chi'}(z)}{\zeta_{\chi'\chi}(z^{-1})}=q^{(g-1)rs}\frac{\LHom(\chi,\chi';z^{-1})}{\LHom(\chi,\chi';qz^{-1})}
$$
and that the sets of zeros $\{u^{\chi,\chi'}_1,\dots,u^{\chi,\chi'}_{2(g-1)rs}\}$ of the polynomials $\zeta_{\chi\chi'}(z)$ satisfy the following relations for all $\chi, \chi' \in \Sigma^{abs}$:
\begin{equation}\label{eqn:functionaleq}
\left\{u^{\chi,\chi'}_1,\dots,u^{\chi,\chi'}_{2(g-1)rs} \right\} = \left\{\frac 1{qu^{\chi',\chi}_1},\dots, \frac 1{qu^{\chi',\chi}_{2(g-1)rs}} \right\}
\end{equation}
Special care must be taken with the formula above when $\chi = \chi'$, in which case the zeroes of the polynomial $\zeta_{\chi\chi}(z) \frac {1-z}{1-zq^{-1}}$ enjoy the symmetry property above.

\medskip

\subsection{} \noindent Let $\underline{\Sigma}^{abs}=\{\underline{\chi}\}$ be a fixed set of representatives  of $\Sigma^{abs}/\mathbb{C}^*$, and fix a corresponding set $\{f_{\underline{\chi}}\}$ of cuspidal eigenforms (choosing different representatives will modify the rational functions defined in the previous Subsection according to \eqref{eqn:rescale}). To this datum, we associate the quiver $Q_X^{abs}$, with vertex set:
$$
I=\bigsqcup_{r = 1}^{\infty} \underline{\Sigma}_r^{abs}
$$
\footnote{Even though the set $I$ is countable, the results in the present paper still hold as stated, because all direct summands of \eqref{eqn:big shuffle} and all relations \eqref{eqn:serre a} only involve finitely many elements of $I$.} and edge set ${E}$ defined as follows~:
\begin{enumerate}

\item[-] if $({\underline{\chi}}, {\underline{\chi}'}) \in \underline{\Sigma}_r^{abs} \times \underline{\Sigma}_s^{abs}$ are distinct, then we draw $(g-1)rs$ arrows going from ${\underline{\chi}}$ to ${\underline{\chi}'}$

\item[-] if ${\underline{\chi}} \in \underline{\Sigma}_r^{abs}$ then there are $(g-1)r^2 + 1$ loops at ${\underline{\chi}}$.
\end{enumerate}
From the definition of $\zeta_{\underline{\chi} \underline{\chi'}}(z)$ and properties ii), iii)  and \eqref{eqn:functionaleq} of Rankin-Selberg L-functions, we see that the function $\zeta_{\underline{\chi} \underline{\chi'}}(z)$ is, up to the (nonzero) constant factor $\theta_{\chi,\chi'}$, exactly of the form \eqref{eqn:def zeta}, upon the specialization of the parameters: 
\begin{equation}
\label{eqn:eq parameters hall}
\left\{t^{\underline{\chi},\underline{\chi}'}_e \right\} \cup \left\{q/t_e^{\underline{\chi}',\underline{\chi}}\right\}=\left\{1/u^{\underline{\chi},\underline{\chi}'}_e \right\}, \qquad (\underline{\chi},\underline{\chi}') \in (\underline{\Sigma}^{abs})^2.
\end{equation}
Let us denote by ${_\BK\CS}$ the associated (small) shuffle algebra, as in \eqref{eqn:restricted wheel}. 

\medskip

\begin{theorem}The assignment $z_{\underline{\chi}}^d \mapsto f_{\underline{\chi},d}$ for $\underline{\chi} \in \underline{\Sigma}^{abs}$ extends to an algebra isomorphism:
$$
{_\BK\CS}  \xrightarrow{\sim}  \mathbf{H}^{+,abs}
$$
\end{theorem}

\medskip

\begin{proof}
This is essentially \cite[Theorem 3.10]{KSV} (the shuffle kernels differ by the factors $\theta_{\underline{\chi},\underline{\chi'}}$, but this does not affect any of the arguments in the present paper).
\end{proof}

Note that \cite{KSV} does not single out absolutely cuspidal eigenforms; however, the kernels for the shuffle algebras considered there are, for not absolutely cuspidal eigenforms, different from the ones considered in the current paper.

\medskip

\noindent In combination with Theorem~\ref{thm:main restricted}, we deduce the following presentation for $\mathbf{H}^{abs}$. Let us set, for $x \in X(\finite)$ and $ \underline{\chi} \in \underline{\Sigma}^{abs}$~:
$$
E_{\underline{\chi}}(z)=\sum_{d \in \BZ} f_{\underline{\chi},d} z^{-d}, \qquad T_x(z)=\sum_{l=1}^{\infty} T_{x,l} z^{-l}
$$

\medskip

\begin{corollary} 
\label{cor:abscusp}

The algebra $\mathbf{H}^{abs}$ is isomorphic to the algebra generated by elements $\kappa$, $\{T_{x,l}\;|\; x \in X(\finite), l \geq 1\}$ and $\{f_{\underline{\chi},d}\;|\; \underline{\chi} \in \underline{\Sigma}^{abs}, d \in \mathbb{Z}\}$ modulo the following set of relations:
\begin{align*} 
&\kappa \text{ and } T_{x,l} \text{ all commute } (x \in X(\finite), l \geq 1) \\[4pt]
&\kappa E_{\underline{\chi}}(z)=q^{(g-1)r}E_{\underline{\chi}}(z)\kappa\ , \qquad (\underline{\chi} \in \Sigma^{abs}_r)\ , \\[4pt]
&[T_x(z),E_{\underline{\chi}}(w)]=\left( \sum_{l = 1}^{\infty} \underline{\chi}(T_{x,l})\left(\frac{w}{z}\right)^l\right) E_{\underline{\chi}}(w)\ , \\[4pt]
&E_\chi(z)E_{\chi'}(w)\zeta_{\chi'\chi} \left( \frac wz \right)=E_{\chi'}(z)E_{\chi}(w)\zeta_{\chi\chi'} \left( \frac zw \right)
\end{align*}
and the collection of cubic relations determined by setting \eqref{eqn:defX restricted} equal to 0. 

\end{corollary}

\medskip

\noindent In particular, all the relations satisfied by Eisenstein series attached to absolutely cuspidal eigenforms are a consequence, in addition to the usual functional equation, of certain \textit{cubic} relations (and the structure of these only depend on the multiplicities of the zeros of the corresponding Rankin-Selberg L-functions).
 
\medskip

\begin{remark}\label{rem:absHall} One might wonder if there is an analogue for $\mathbf{H}_X$ of the \emph{generic} spherical Hall algebra. In the context of a quiver $Q$, such an analogue takes the form of a quantized \emph{graded} Borcherds algebra, whose Cartan datum encodes the dimensions of the spaces of cuspidal functions in $\mathbf{H}_Q$, see \cite{BoS}. By a recent theorem of H. Yu, \cite{HongjieYu}, the dimensions of the spaces of \emph{absolutely} cuspidal functions in $\mathbf{H}_X$ are given by some universal polynomials in the Weil numbers of $X$ (depending on the rank of the sheaves considered). This suggests that it is $\mathbf{H}_X^{abs}$ rather than $\mathbf{H}_X$ which admits a natural generic form, and that the classical limit of such a generic form of $\mathbf{H}_X^{abs}$ would be a Lie algebra in the category of $GSp_{2g}(\mathbb{C})$-modules. However, we do not know at the moment how to encode the zeros of the various Rankin-Selberg L-functions (which account for the ``quantum" parameters).
\end{remark}

\bigskip

\newcommand{\etalchar}[1]{$^{#1}$}
\providecommand{\bysame}{\leavevmode\hbox to3em{\hrulefill}\thinspace}
\providecommand{\MR}{\relax\ifhmode\unskip\space\fi MR }
\providecommand{\MRhref}[2]{%
	\href{http://www.ams.org/mathscinet-getitem?mr=#1}{#2}
}
\providecommand{\href}[2]{#2}

\end{document}